\documentclass[11pt,a4paper]{article}

\usepackage{amsfonts}
\usepackage{amsmath}
\usepackage{graphicx}
\usepackage{amssymb}

\usepackage{color}
\usepackage{algorithm}
\usepackage{algorithmic}

\usepackage{enumitem}
\usepackage{multirow}

\numberwithin{equation}{subsection}

\newtheorem{thm}{Theorem}[section]

\newcommand{\eps}{\varepsilon}

\title{Overlapping non Matching Meshes Domain Decomposition Method in Isogeometric Analysis}
\author{Michel Bercovier \thanks{Corresponding author. E-Mail :\texttt {berco@cs.huji.ac.il};}   and Ilya Soloveichik \\  \small{The Rachel and Selim Benin School of Computer Science and Engineering} ,\\ \small{Hebrew University of Jerusalem,Israel.}  }

\begin{document}
\parskip=.050in

%
\renewcommand{\thepage}{\roman{page}}
\maketitle
\addcontentsline{toc}{section}{Abstract}

\noindent

\begin{abstract}
One of the important aspects of IsoGeometric Analysis (IGA) is its link to Computer Aided Design (CAD) geometry methods.
Two of IGA's  major challenges are geometries made of boolean assemblies  such as in
Constructive Solid Geometry (CSG), and local refinements.  In the present work we propose to
 apply the Additive  Schwarz Domain Decomposition Method  (ASDDM) on overlapping  domains to actually solve Partial Differential Equations (PDEs), on assemblies as well as a means for local refinement.
 
As a first step we study a collection of simple two domains problems, using the software package GeoPDEs,  all of them converged rapidly, even on very distorted domains.
For local refinement we implement a "Chimera" type of zooming together with ASDDM iterations, and here again we get very attractive results.
A last effect is that the additive  methods are naturally parallelizable and thus can extend to large
3D problems. We give some examples of computations on  multi-patch bodies.

This new implementation of the DD methods brings many interesting  problems such as: the definition of boundary conditions on trimmed patches, choice of blending operators on the intersections ,  optimization  of the overlap size, choice of preconditionning to guarantee convergence ( as we do not have a maximum principle here.)
\end{abstract}
\eject

\vspace{-0.3in}
\tableofcontents

\eject
\setcounter{page}{1}
\renewcommand{\thepage}{\arabic{page}}

\section{Motivation: Using IsoGeometry and Domain Decomposition Methods Together}

\begin{figure}
\includegraphics[width=0.6\textwidth]{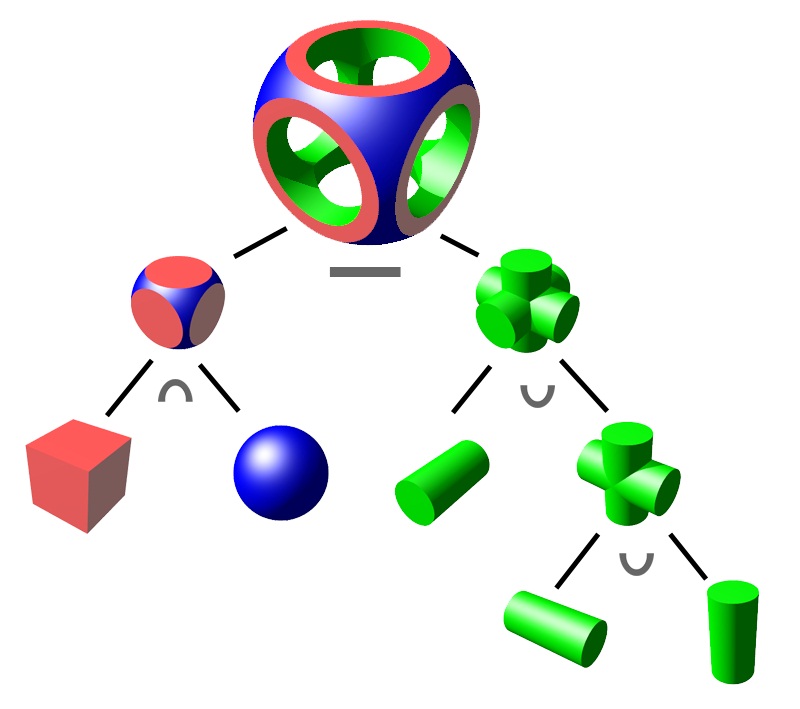}
\caption{CSG example of Boolean operations.}
\label{fig:csg}
\end{figure}

One of the main purposes of Constructive Solid Geometry is to provide CAD with the tools allowing the engineer to perform Boolean operations on the bodies he works with. There are very many engineering forms which are obtained as unions, intersections, subtractions etc. of basic geometrical figures, see Figure \ref{fig:csg} for example. Even quite complicated bodies can often be represented as results of such set operations applied to some primitive geometrical forms like cubes, balls, cylinders etc. \cite{hoffman}. To solve Partial differential equations on domains defined by such solids , it seems natural to find methods that take advantage of the definitions of such solids.

IsoGeometric Analysis (IGA) ~\cite{hughes} is a relatively new numerical method of solving PDEs, derived from the FEM. The main features of IGA, which are very important in engineering and which FEM lacks, are: (1) exact representation of the geometry and (2) the same family of basis functions that are used in the construction of the geometry are used in the solution of PDEs, thus avoiding the cumbersome mesh generation step. 

\begin{figure}
\includegraphics[width=0.7\textwidth]{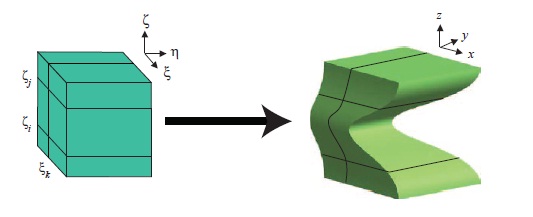}
\caption{An example of IsoGeometric mapping.}
\label{fig:intro_pushforward}
\end{figure}

The main concept of IGA can be described in the following way. Assume that we work in 3-dimensional environment and consider a cube, which constitutes the reference space. The geometrical domain is build as a continuous deformation of the reference cube, as shown in Figure ~\ref{fig:intro_pushforward}, and represents the exact geometrical form. Given a finite space of basis functions on the reference domain we apply the geometrical mapping to them to generate the finite space of solutions, as the space spanned by their images.

Since a single patch is limited, one would like to extend  IGA to complex domainss, and based on CSG and consider a collection of subdomains . Domain Decomposition (DD) methods are then natural candidates for the solution of PDEs over complex domains.  More for  IGA  the problem of local zooming becomes more involved since the refinement need to be introduced at the level of the parametric space, where the rigid structure of tensor product must be respected.( In view of this restriction different methods of local zooming were proposed for IGA. Popular ways of doing that is the T-spline technique ~\cite{tspline} or the Truncated Hierarchical B Splines (THB) one, ~\cite{thb}. This allows one to locally refine the mesh, but still the tensor product structure of the parametric space leads to non-local changes in the mesh ~\cite{dorf}, implying that many additional degrees of freedom are introduced around the areas being magnified, and increasing the computational costs of  the solution. ).  Domain Decomposition (DD) methods, can also bring a solution here by using local zooming or Chimera type methods introduced  in  Dougherty at all ~\cite{steger}.

Domain Decomposition partitions the original domain into a set of overlapping( or just glued together) connected subdomains with smooth boundary and then solve the PDEs on each of them iteratively, see ~\cite{widlund} for a detailed exposition. In most implementations the meshes of contiguous domains match or are linked by a mortar element.
It has been studied for IGA in ~\cite{cho}  and ~\cite{DD-IGA}.

In the present work we introduce the definition and properties  of DD methods for non-matching meshes over \emph {overlapping} subdomain meshes as well as  applications to problems of local refinement (zooming) . For zooming the meshes of the newly constructed subdomains are independently defined  and \emph{non-matching }, but the geometry of the whole domain is  unchanged. This way of solving the local zooming problem allows us to introduce  new degrees of freedom only in the region under investigation, at the same time we use the IGA framework with all its advantages.

\section{Overview of this work}

In the next part we recall some   definitions and notations related to B-splines, which are the cornerstones of all the IGA technique. We consider the refinement procedures such as knot insertion and degree elevation, which allow us to introduce local refinement and zooming to the IGA. 

As a model example we state  the classical Poisson equation in its continuous form. We then proceed to the weak form of this problem, which gives us quite a natural way of solving the equation numerically,  the Galerkin method, used in FEM and IGA. Then the main ideas of IGA are explained with all the computations needed to perform the actual calculations of the approximate solution. Imposing non homogeneous Dirichlet boundary conditions is not straight forward in IGA and some special techniques should be used to impose them. We describe two different methods for imposing essential boundary conditions. 

The next part introduces the family of Schwarz Domain Decomposition methods and gives their description. We choose to implement the Additive Schwarz Domain Decomposition Method (ASDDM), that is equivalent to a block Jacobi iterative process since it may be easily parallelized and we detail its continuous and discretized weak forms.

In the sixth part the Discretized Schwarz Domain Decomposition Method is described. The construction of exact and approximation projection operators for the imposition of Dirichlet boundary conditions is treated.

The next part  illustrates the Schwarz Algorithm in the one-dimensional case. We give the detailed overview of the matrix operator forms obtained from the ASDDM.

We then give  numerous  numerical results. Based on the GeoPDEs software ~\cite{geopdes}, we have developed a code that implements the ASDDM using IGA and applied it to solve different one-, two-, and three-dimensional problems. In this section all the convergence results are provided. We also describe how the code was parallelized and show three-dimensional examples, computed by it. 

We conclude by discussing further directions for this  research and open problems.

\section{B-splines }
\subsection{Introduction}
The most popular geometrical forms being used in CAGD are based on B-splines. Their definition provides us with very important geometrical properties. In particular, the degree of a B-spline curve is strictly related to the number of control points and they give us local control over the form of the curve. 
B-splines are not general enough to describe exactly even very simple shapes such as cylinders, balls, circles, etc.
The notion of Non -
Uniform Rational B-Splines (NURBS) extends B-splines and allows one to exactly represent a wide array of objects that cannot be exactly represented by polynomials,
many of which are ubiquitous in engineering design.  In the present work we will consider only B-Splines .We will now provide the exact definitions. All the details can be found in ~\cite{peigl}.

\subsection{B-spline basis functions}

Let $\Xi$ be a set of $m$ non-decreasing numbers, $\xi_1 \leq \xi_1 \leq
\xi_2 \leq \dots \leq \xi_m$. The $\xi_i$'s are called $\textit{knots}$, the set
$\Xi$ is the $\textit{knot}$ $\textit{vector}$, and the half-open interval $[\xi_i, \xi_{i+1})$ is the
$i$-th $\textit{knot}$ $\textit{span}$. If a knot $\xi_i$ appears $k$ times (i.e.,
$\xi_i = \xi_{i+1} = \dots = \xi_{i+k-1}$), where $k > 1$, $\xi_i$ is a
multiple knot of multiplicity $k$ and the corresponding knot span does not exist. Otherwise,
if $\xi_i$ appears only once, it is a simple knot. If the knots are
equally spaced, the knot vector or the knot sequence is said to be uniform;
otherwise, it is non-uniform.

\begin{figure}

\includegraphics[width=0.8\textwidth]{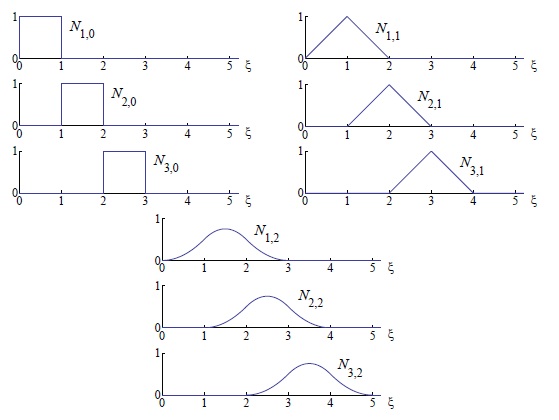}
\caption{B-spline basis functions of order 0, 1 and 2 for uniform knot vector $\Xi = [0,1,2,3,\dots]$.}
\label{fig:b_spline_basis}
\end{figure}

To define B-spline basis functions, we need one more parameter, the
degree of these basis functions, $p$. The $i$-th B-spline basis
function of degree $p$, written as $N_{i,p}(\xi)$, is defined
recursively as follows:

\begin{equation}
  N_{i,0}=\left\{ 
  \begin{array}{l l}
    1 & \quad \text{if $\xi_i \leq \xi < \xi_{i+1}$,}\\
    0 & \quad \text{otherwise.}\\
  \end{array} \right.
  \label{eq:bspline_def1}
\end{equation}

\begin{equation}
N_{i,p}(\xi)=\frac{\xi-\xi_i}{\xi_{i+p}-\xi_i}N_{i,p-1}(\xi)+\frac{\xi_{i+p+1}-\xi}{\xi_{i+p+1}-\xi_{i+1}}N_{i+1,p-1}(\xi).
\label{eq:bspline_def2}
\end{equation}

\begin{figure}
\centering
\includegraphics[width=0.7\textwidth]{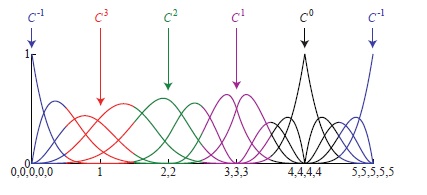}
\caption{Quartic $(p=4)$ basis functions for an open, non-uniform knot vector $\Xi = [0,0,0,0,0,1,2,2,3,3,3,4,4,4,4,5,5,5,5,5]$. The continuity across an interior element boundary is a direct result of the polynomial order and the multiplicity of the corresponding knot value.}
\label{fig:b-spline_basis}
\end{figure}

If the degree is zero (i.e., $p = 0$), these basis
functions are all indicators of the the intervals they are defined on. 

In general, basis functions of order $p$ have $p-k$ continuous derivatives across knot $\xi_i$, where $k$ is the multiplicity of the value $\xi_i$ as before. When the multiplicity of a knot value is exactly $p$, the basis is interpolatory at that knot. When the multiplicity is $p+1$, the basis becomes discontinuous; the patch boundaries are formed in this way.

\subsection{B-spline curves}
In order to define a B-spline we have to provide the following information: a set of $n$ control points $\{P_i\}_{i=1}^n$, a knot vector $\Xi = \{\xi_i\}_{i=1}^m$ of $m$ knots, and a degree $p$, such that $n$, $m$ and $p$ must satisfy $m = n + p
+ 1$. 

Given this information the B-spline curve of
degree $p$ defined by these control points and knot vector $\Xi$ is
$$
\mathbf{C}(\xi)=\sum_{i=1}^n N_{i,p}(\xi)\mathbf{P}_i,
$$
where $N_{i,p}(u)$'s are the B-spline basis functions of degree $p$.

The point on the curve that corresponds to a knot $\xi_i$, $C(\xi_i)$, is referred to as
a $\textit{knot}$ $\textit{point}$. Hence, the knot points divide a B-spline curve into
curve segments, each of which is defined on a knot span.

\begin{figure}
\centering
\includegraphics[width=0.6\textwidth]{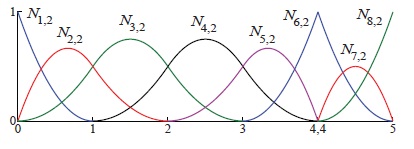}
\caption{Quadratic basis functions for open, non-uniform knot vector $\Xi = \{0, 0, 0, 1, 2, 3, 4,
4, 5, 5, 5 \}$}
\label{fig:basis_curve}
\end{figure}

\begin{figure}
\centering
\includegraphics[width=0.4\textwidth]{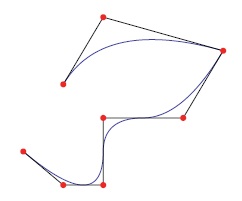}
\caption{An example of a B-spline curve and its control points}
\label{fig:curve}
\end{figure}

If the knot vector does not have any particular structure, the
generated curve will not touch the first and last legs of the
control polygon. This type of
B-spline curves is called open B-spline curves. In order to make the curve interpolatory the first knot and the last knot must be of
multiplicity $p+1$. By repeating the first and the last knots $p+1$ times the curve starts at $P_0$ and ends at $P_n$. A knot vector of such form is called an open knot vector. Here we consider all the knot vectors to be open if we don\rq{}t state otherwise (e.g., it is not the case when periodic NURBS are discussed).

In Figure ~\ref{fig:curve} you can see an example of a B-spline curve. We chose a non-uniform knot vector $\Xi = \{0, 0, 0, 1, 2, 3, 4,
4, 5, 5, 5 \}$ and the control points as shown in the figure. The corresponding parametric space and the B-spline basis functions are shown in Figure ~\ref{fig:basis_curve}.

\subsection{Derivatives of B-spline }
Since the B-spline basis functions are obtained by  recursion from the lower degree B-spline basis functions and the dependence is linear, we can obtain the recursive formulas for their derivatives in a similar way and represent them as a linear combination of B-spline basis functions of lower degrees, see ~\cite{peigl}.

The derivative of a basis function is given by
\begin{equation}
N\rq{}_{i,p}(\xi)=\frac{p}{\xi_{i+p}-\xi_i}N_{i,p-1}(\xi)-\frac{p}{\xi_{i+p+1}-\xi_{i+1}}N_{i+1,p-1}(\xi).
\label{eq:b-spline_der}
\end{equation}
The proof of this recurrence is obtained by induction on $p$.

It is a little more complicated task to calculate the NURBS derivatives since they are fractions of linear combinations of the B-spline basis functions, but still the recurrence relation exists and may be found in the Appendix.

\subsection{Multi-dimensional B-splines }
All the definitions given above refer to curves. One can easily extend these definitions to surfaces and bodies by taking a tensor product of an appropriate number of knot vectors and defining the basis functions as the products of the corresponding one-dimensional basis functions. More exactly, given a control net $\{\mathbf{P}_{i,j}\}_{i,j=1}^{n,m}$, polynomial orders $p$ and $q$, and knot vectors $\Xi=\{\xi_i\}_{i=1}^{n+p+1}$ and $\Upsilon=\{\upsilon_i\}_{i=1}^{m+q+1}$, a tensor product B-spline surface is defined by

\begin{equation}
\mathbf{S}(\xi,\upsilon) = \sum_{i=1}^{n} \sum_{j=1}^{m} N_{i,p}(\xi) M_{j,q}(\upsilon) \mathbf{P}_{i,j},
\end{equation}
where $N_{i,p}(\xi)$ and $M_{j,q}(\upsilon)$ are the univariate B-spline basis functions corresponding to knot vectors $\Xi$ and $\Upsilon$, respectively.

Many of the properties of a B-spline surface are the result of its
tensor product nature. The basis functions are point-wise nonnegative, and form
a partition of unity as for any $(\xi,\upsilon) \in
[\xi_1,\xi_{n+p+1}]\times [\upsilon_1,\upsilon_{m+q+1}],$
$$
\sum_{i=1,j=1}^{n,m} N_{i,p}(\xi) M_{j,q}(\upsilon) = \left(
\sum_{i=1}^{n} N_{i,p}(\xi) \right) \left( \sum_{j=1}^{m}
M_{j,q}(\upsilon) \right) = 1.
$$

The local support of the basis functions also follows directly from
the one-dimensional functions that form them. The support of a given
bivariate function $\widetilde{N}_{i,j;p,q}(\xi,\upsilon) =
N_{i,p}(\xi) M_{j,q}(\upsilon)$ is exactly
$[\xi_i,\xi_{i+p+1}]\times [\upsilon_j,\upsilon_{j+q+1}]$

Tensor product B-spline solids are defined in analogous fashion to
B-spline surfaces. Given a control lattice $\{\mathbf{P}_{i,j,k}\}_{i,j,k=1}^{n,m,l}$, 
polynomial orders $p,q$
and $r$, and knot vectors $\Xi = [\xi_1, \xi_2,\dots,\xi_{n+p+1}],
\Upsilon = [\upsilon_1, \upsilon_2,\dots,\upsilon_{m+q+1}],$ and
$\Sigma = [\sigma_1, \sigma_2,\dots,\sigma_{l+r+1}],$ a B-spline
solid is defined by
$$
\mathbf{S}(\xi,\upsilon) = \sum_{i,j,k=1}^{n,m,l} N_{i,p}(\xi)
M_{j,q}(\upsilon) L_{k,r}(\sigma) \mathbf{P}_{i,j,k}.
$$

The properties of a B-spline solid are trivariate generalizations of
those for B-spline surfaces.

\begin{figure}
\centering
\includegraphics[width=0.7\textwidth]{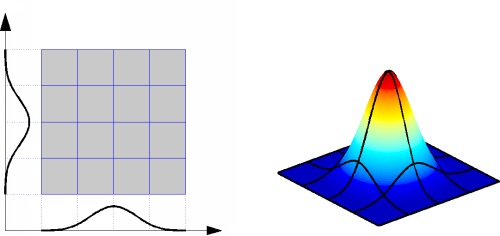}
\caption{An example of a two-dimensional B-spline basis function, which is a tensor product of two one-dimensional B-splines}
\label{fig:2d_basis_1}
\end{figure}

One of the most important properties of B-splines is the number of ways in which the basis may be enriched while leaving the underlying geometry and its parameterization unaffected. We consider here two different ways of doing such enrichment: knot insertion and degree elevation. We also want to mention that not only do we have control over the element size and the order of the basis, but we can also control the continuity of the basis as well. We refer to ~\cite{peigl} for the details and ~\cite{hughes} for the applications to IGA.

\begin{figure}
\centering
\includegraphics[width=0.7\textwidth]{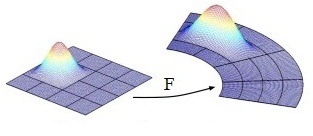}
\caption{$\mathbf{F}$-mapping and a basis function transformation}
\label{fig:F_map_basis_2}
\end{figure}

\begin{figure}
\centering
\includegraphics[width=0.65\textwidth]{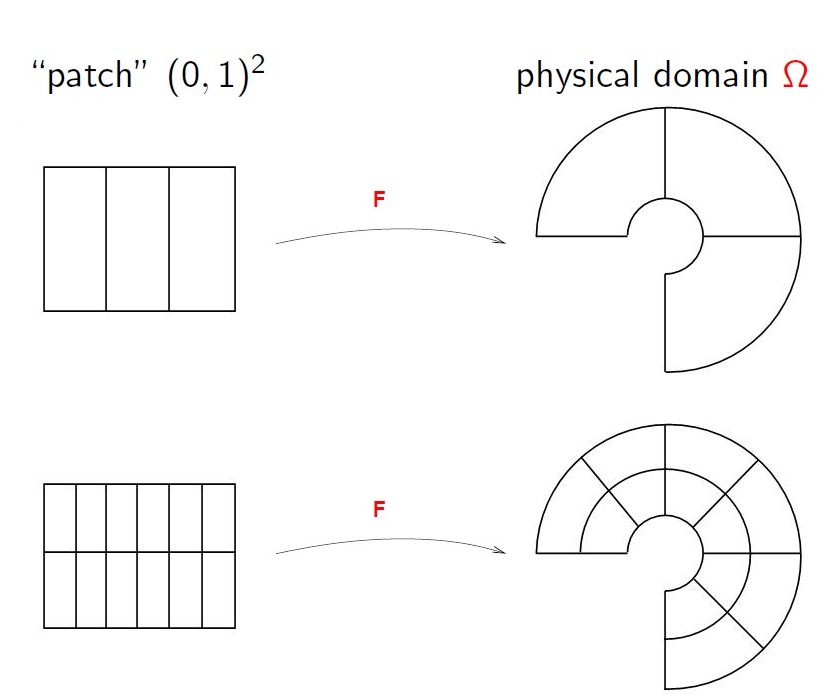}
\caption{Example of a $\mathbf{F}$-mapping and refinement}
\label{fig:Fmapping_ref}
\end{figure}

\subsection{The IsoGeometric Concept}
The fundamental idea \ behind IGA is that the basis used to
exactly model the geometry will also serve as the basis for the
solution space of the numerical method. This notion of using the
same basis for geometry and analysis is called the
$\textbf{Isoparametric  concept}$, and it is quite common in
classical finite element analysis. The fundamental difference
between the new concept of Isogeometric Analysis and the classic concept
of isoparametric finite element analysis is that, in classical FEA,
the basis chosen to \textit{approximate} the unknown solution fields is also
used to define the domain . Isogeometric Analysis turns this
idea around and selects a basis capable of exactly representing the
known geometry and uses it as a basis for the field we wish to
approximate. In a sense, we are reversing the isoparametric arrow
such that it points from geometry toward the solution space, rather
than vice versa.

One of the main concepts of IGA is the geometrical mapping $\mathbf{F}$. It acts on the parametric domain $\widehat{\Omega}$:
\begin{equation}
\mathbf{F}: \widehat{\Omega}\rightarrow\Omega,
\label{eq:parameterization}
\end{equation}
defining  the geometrical domain in the physical space. In fact it acts on the basis functions of the parametric space defining both the geometry and the space of approximation functions used for approximating the solutions of PDEs. Figure \ref{fig:F_map_basis_2} shows how the shape is constructed and the transformation of a basis function under the $\mathbf{F}$-mapping. In the present work the basis functions will always be NURBS or B-splines.

\subsection{Refinement and IGA} One of the important features of IGA is that when refining the parametric domain the geometrical mapping remains the same resulting in preserving the geometry while refining the solution space. Thus, IGA allows us to make local refinements without changing the geometrical shape. An example of a two-dimensional  knot-insertion is shown in Figure \ref{fig:Fmapping_ref}.

One of the drawbacks of IGA is the nonlocal nature of the refinement due to the tensor product structure of the B-splines. The difficulties arise in large scale problems since the number of degrees of freedom grows very fast.  Different algorithms, such as T-splines ~\cite{dorf}, were introduced into IGA in order to avoid these problems. All of them still makes the refinement non-local since some tensor product structure needs to be preserved.

This led us to think about different methods of local mesh refinement, which would allow us to keep the computational costs low. We propose to use the Domain Decomposition Method. Its main difference from the refinement techniques is that the new refined mesh is an absolutely separated domain and does not match the coarse one, but now we have to perform an iterative procedure in order to obtain the solution.

\subsection{The Alternating Schwarz Method}
The Alternating Schwarz method belongs to the class of domain decomposition methods for partial differential equations. Domain decomposition methods can be regarded as divide and conquer algorithms. The main idea is the following: given an open domain (we will assume it to possess some additional properties in the sequel, like being smooth and connected) $\Omega$ we partition it into a number of subdomains $\{\Omega_i\}_{i=1}^n$, satisfying some conditions, and such that $\Omega = \bigcup_{i=1}^n \Omega_i$. The original problem then can be reformulated as a family of subproblems of reduced size and complexity defined on the subdomains.

One of the major advantages of domain decomposition methods is the natural parallelism in solving the subproblems defined on the subdomains. With the upcoming of parallel computer architectures, domain decomposition methods have become very popular during the last two decades. However the origin goes back to 1870. In ~\cite{Sc70}, H. A. Schwarz introduced an algorithm to prove the existence of harmonic functions on irregularly shaped domains. Today this algorithm is known as  the Alternating Schwarz Method.

In the next section we will define the problem we are going to solve and the main tool we use to solve it:  the IGA framework.

\section{Notations and the Problem formulation}

Given a bounded connected open domain $\Omega$ with Lipshitz continuous $\partial{\Omega}$, let us consider the Poisson equation
\begin{eqnarray}
- \text{\textbf{div}}{(k(x) \text{\textbf{grad}}{u})} = f & \quad \text{in } \Omega, \label{eq:poisson}\\
k(x) \frac{\partial{u}}{\partial{\mathbf{n}}} = h & \quad \text{ on } \Gamma_N, \label{eq:neum_cond}\\
u = g & \quad \text{ on }\Gamma_D, \label{eq:dir_cond}
\end{eqnarray} 
where $\overline{\Gamma_D \bigcup \Gamma_N} =
\partial \Omega$, $\Gamma_D \text{ and } \Gamma_N$ are disjoint, and $\mathbf{n}$ is the unit outward normal
vector on $\partial \Omega$. The functions $f \in \mathcal{H}^{-1}(\Omega), g \in \mathcal{L}^{2}(\Gamma_D)$ representing the Dirichlet boundary conditions and $h \in \mathcal{L}^2(\Gamma_N)$ representing the Neumann boundary conditions, are all given.

For a sufficiently smooth domain, a unique solution $u$ satisfying $\eqref{eq:poisson} - \eqref{eq:dir_cond}$ is known to
exist under some standard conditions.

\subsection{Weak form of the problem}

The numerical technique we are going to take advantage of begins by
defining a weak, or variational, counterpart of $\eqref{eq:poisson} - \eqref{eq:dir_cond}$. To do so, we
need to characterize two classes of functions. The first is to be
composed of candidate, or trial functions. From the outset, these
functions will be required to satisfy the Dirichlet boundary
conditions of $\eqref{eq:dir_cond}$.

We may now define the collection of $\textbf{trial functions}$, denoted by
$\mathcal{S}$ as
$$
\mathcal{S} = \{u:u \in \mathcal{H}^1(\Omega), u|_{\Gamma_D} = g\}.
$$

The second collection of functions in which we are interested is
called the $\textbf{test}$ or $\textbf{weighting functions}$. This collection is very similar to
the trial functions, except that we have the homogeneous counterpart
of the Dirichlet boundary conditions:
\begin{equation}
\mathcal{V} = \{w:w \in \mathcal{H}^1(\Omega), w|_{\Gamma_D} = 0\},
\label{eq:hom_sp}
\end{equation}

which is a Hilbert space with respect to the associated energy norm $ ||u|| =a(u,u)^{1/2}$, that follows from the coercivity of the bilinear form $a(\cdot,\cdot)$, as shown below.

We may now obtain a variational statement of our boundary value
problem by multiplying the equation $\eqref{eq:poisson}$ by an arbitrary
test function $w \in \mathcal{V}$ and integrating by parts,
incorporating $\eqref{eq:neum_cond}$ as needed. 

Given $f, g$ and $h$ find $u \in \mathcal{S}$
such that for all $w \in \mathcal{V}$
\begin{equation}
\int_{\Omega} k(x) \textbf{grad} w \cdot \textbf{grad} u \textbf{d}
\Omega = \int_{\Omega}
w f \textbf{d} \Omega + \int_{\Gamma_N} w h \textbf{d} \Gamma.
\label{eq:gen_weak}
\end{equation}

Define a bilinear form $a(\cdot,\cdot): \mathcal{H}^1(\Omega) \times \mathcal{H}^1(\Omega) \rightarrow \mathbb{R}$ and a functional $L: \mathcal{H}^1(\Omega) \rightarrow \mathbb{R}$ as:
\begin{equation}
a(w,u) = \int_{\Omega} k(x) \textbf{grad} w \cdot \textbf{grad} u
\textbf{d} \Omega,
\end{equation}
and
\begin{equation}
L(w) = \int_{\Omega}  w f \textbf{d} \Omega + \int_{\Gamma_N} w h
\textbf{d} \Gamma.
\end{equation}
Now the weak form reads as:
\begin{equation}
a(w,u) = L(w),
\label{eq:var_form}
\end{equation}

We shall note here, that the bilinear functional $a(\cdot,\cdot)$ is
symmetric and positive definite for any $k(x)>0$ when considered over the space $\mathcal{V}$. Indeed, as shown below it is coercive, which results in the symmetry and positive definiteness of the stiffness matrix we are going to define below . One will find in  (CIARLET)   that the above problem ~\ref{eq:var_form}  has a unique solution in  $\mathcal{H}^1(\Omega)$.

\subsection{Galerkin's Method}

In Galerkin's method we define finite-dimensional
approximations of the spaces $\mathcal{S}$ and $\mathcal{V}$, denoted
$\mathcal{S}_h$ and $\mathcal{V}_h$ such that $\mathcal{S}_h \subset \mathcal{S},$ 
$\mathcal{V}_h \subset \mathcal{V}$. These subspaces will be associated with subsets of the
space spanned by the isogeometric basis.

As we have just mentioned we can further characterize $\mathcal{S}_h$ by recognizing that if
we have a given function $g_h \in \mathcal{S}_h$ such that
$g_h|_{\Gamma_D} = g$, then for every $u_h \in \mathcal{S}_h$ there
exists a unique $v_h \in \mathcal{V}_h$ such that
$$
u_h = v_h + g_h.
$$

This clearly will not be possible for an arbitrary function $g$. The
problem of approximation of the function $g$ by a function $g_h \in
\mathcal{S}_h$ will be treated in detail below. So let us assume at
present that such a function $g_h$ exists.

As in  ~\cite{hughes} we can write the Galerkin\rq{}s form of the problem as: given $g_h$ and $h$, find $u_h =
v_h + g_h$, where $v_h \in \mathcal{V}_h$, such that for all $w_h
\in \mathcal{V}_h$
\begin{equation}
a(w_h,u_h) = L(w_h),
\label{eq:discr_var_form}
\end{equation}
or
\begin{equation}
a(w_h,v_h) = L(w_h) - a(w_h,g_h).
\end{equation}

We will use this form for the variational formulation for the Schwarz Method.

In general the weighting space $\mathcal{V}_h$ of test functions may be different
than the space $\widetilde{\mathcal{V}}_h$ of trial functions $v_h$, that is $v_h \in \mathcal{V}_h$ but $w_h \in
\widetilde{\mathcal{V}}_h \neq \mathcal{V}_h$, but in the following
we will assume they are the same.

\section{General IGA}
\label{sec:General_IGA}

The goal of IGA, as it is also for FEM, is the numerical
approximation of the solution of partial differential equations
(PDEs) \cite{geopdes}. In both approaches the PDEs are numerically solved using the
Galerkin procedure: the equations are written in there equivalent
variational formulations, and a solution is sought in a finite
dimensional space with good approximation properties. The main
difference between the two methodologies is that in FEM the basis
functions and the computational geometry (i.e., the mesh) are
defined element by element, using piecewise polynomial approximations, whereas in IGA the
computational geometry and the space are defined \textit{exactly} from the (input data) information and
the basis functions (e.g., NURBS, T-splines or generalized
B-splines) given by CAD.

Another important difference is that the unknowns to be computed, coefficients that define
the solution as a linear combination of the basis functions (degrees of freedom), are not nodal values, but \lq\lq{}control values\rq\rq{} .

\begin{figure}
\centering
\includegraphics[width=0.7\textwidth]{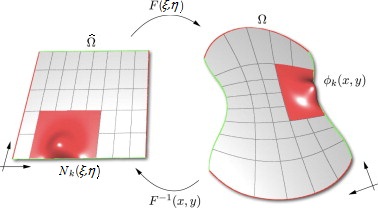}
\caption{Example of a $\mathbf{F}$-mapping and the basis function transformation}
\label{fig:f_map_basis}
\end{figure}

Let us consider a two-dimensional case where we assume that the
physical domain $\Omega \subset \mathbb{R}^2$ is open, bounded and
piecewise smooth. We also assume that such a domain can be exactly
described through a parametrization of the form
\begin{equation}
\mathbf{F}: \widehat{\Omega}\rightarrow\Omega,
\label{eq:parametrization}
\end{equation}
where $\widehat{\Omega}$ is some parametric domain (e.g., the unit
square), and the value of the parameterization can be computed via the isogeometric mapping $\mathbf{F}$, which is
assumed to be smooth with piecewise smooth inverse, see Figure \ref{fig:f_map_basis}.

We now construct the finite-dimensional spaces $\mathcal{S}_h$ and $\mathcal{V}_h$ as spans over their basis functions:
$$
\mathcal{S}_h := \langle \phi_1,
\phi_2,...,\phi_{N_h},\phi_{N_h+1},...,\phi_{N_h+N_h^b} \rangle
\subset \mathcal{S},
$$
and
$$
\mathcal{V}_h := \langle \phi_1,
\phi_2,...,\phi_{N_h} \rangle
\subset \mathcal{V},
$$
where we first number all the basis function which vanish on the boundary $\Gamma_D$: $\phi_1,
\phi_2,...,\phi_{N_h}$, and then the other basis functions $\phi_{N_h+1},...,\phi_{N_h+N_h^b}$, which do not vanish on this boundary.

The IGA prescribes us the way of building the basis functions for the spaces $\mathcal{S}_h$ and $\mathcal{V}_h$. Let $\{
\widehat{\phi}_j\}_{j \in \mathcal{I} \cup \mathcal{B}}$ be a basis for
$\widehat{\mathcal{S}}_h$, with $\mathcal{I} \cup \mathcal{B}$ a proper set of indices; we partition the indices in the same way as before: $\mathcal{I}$ corresponds to the basis functions that vanish on the boundary $\widehat{\Gamma}_D$ and $\mathcal{B}$ corresponds to those which do not vanish on $\widehat{\Gamma}_D$.
With the assumption made on $\mathbf{F}$, the set
$\{\widehat{\phi}_j \circ \mathbf{F}^{-1} \}_{j \in \mathcal{I} \cup \mathcal{B}}\equiv\{\phi_j\}_{j \in \mathcal{I} \cup \mathcal{B}}$ is a basis for
$\mathcal{S}_h$ and $\{ \widehat{\phi}_j \circ \mathbf{F}^{-1} \}_{j \in \mathcal{I}}\equiv\{\phi_j\}_{j \in \mathcal{I}}$ is a basis for $\mathcal{V}_h$. 

In IGA, as introduced in \cite{hughes}, the spaces $\mathcal{S}_h$ and $\mathcal{V}_h$ are formed by transformation of 
B-spline functions. We define these space in
the following general way:

\begin{equation}
\mathcal{S}_h := \{ u_h \in \mathcal{S} : \widehat{u}_h = u_h \circ \mathbf{F} \in \widehat{\mathcal{S}}_h\}
\equiv \{ u_h \in \mathcal{S} : u_h = \widehat{u}_h \circ \mathbf{F}^{-1}, \widehat{u}_h
\in \widehat{\mathcal{S}}_h \},
\end{equation}
\begin{equation}
\mathcal{V}_h := \{ v_h \in \mathcal{V} : \widehat{v}_h = v_h \circ \mathbf{F} \in \widehat{\mathcal{V}}_h\}
\equiv \{ v_h \in \mathcal{V} : v_h = \widehat{v}_h \circ \mathbf{F}^{-1}, \widehat{v}_h
\in \widehat{\mathcal{V}}_h \},
\end{equation}
where $\iota$ is a proper pull-back, defined from the
parametrization $\eqref{eq:parameterization}$, and $\widehat{\mathcal{S}}_h$ and $\widehat{\mathcal{V}}_h$ are finite spaces generated by the basis functions
defined in the parametric domain $\widehat{\Omega}$.

Hence, the discrete solution of the problem can be
written as

\begin{equation}
u_h = \sum_{j \in \mathcal{I} \cup \mathcal{B}} \alpha_j \phi_j = \sum_{j \in \mathcal{I} \cup \mathcal{B}} \alpha_j \widehat{\phi}_j \circ \mathbf{F}^{-1}.
\label{eq:nonhom_discr_sol}
\end{equation}

We assume that for the function $g$ representing the Dirichlet boundary conditions, there exists a function $g_h \in \mathcal{S}_h$ such that $g_h|_{\Gamma_D} = g$, meaning that there exist $\{ \gamma_j \}_{j \in\mathcal{B}}$ such that 
$g = \sum_{j \in\mathcal{B}} \gamma_j \phi_j$. The more general case of $g \in \mathcal{S}$ will be treated below.

The solution now read as
$$
u_h = \sum_{j \in\mathcal{B}} \gamma_j \phi_j + \sum_{j \in \mathcal{I}} \alpha_j \phi_j.
$$
Thus, instead of considering the function $u_h \in \mathcal{S}_h$ we can treat without loss of generality the function $u_h - \sum_{j \in\mathcal{B}} \gamma_j \phi_j \in \mathcal{V}_h$, reducing the problem to the homogeneous Dirichlet boundary conditions case. So we can write $$
u_h = \sum_{j \in \mathcal{I}} \alpha_j \phi_j.
$$

Substituting this expression into $\eqref{eq:discr_var_form}$, and testing against every
test function $\phi_i \in \mathcal{V}_h$, we obtain a linear system of
equations where the coefficients $\alpha_j$ are the unknowns, and
the entries of the matrix and the right-hand side are
$a(\phi_i,\phi_j)$ and $L(\phi_i)$, respectively. 
\begin{equation}
a({\sum_{j \in \mathcal{I}} \alpha_j \phi_j},\phi_i) = L(\phi_i), \forall i \in \mathcal{I}
\label{eq:weak_discr_problem}
\end{equation}

These terms have to be computed using suitable quadrature rules for
numerical integration \cite{hughes}.

\section{Model Problem: Poisson Equation with homogeneous Dirichlet BC}

We now specialize the general framework of the previous section to the
particular case of the Poisson's problem $\eqref{eq:poisson} - \eqref{eq:dir_cond}$, defined in a physical
domain $\Omega$ described and discretized with B-Splines. This constitutes the model problem on which we show in
detail the IGA method.

Let us assume that the computational domain is constructed as a
singleB-Spline patch, such that the parameterization $\mathbf{F}$ is
given as above, see Figure \ref{fig:f_map_basis}. The assumptions on $\mathbf{F}$ of the previous section are
here supposed to be valid. For the sake of simplicity, homogeneous Dirichlet boundary conditions are assumed.

The linear equations from \ref{eq:weak_discr_problem} can be rewriten as:
\begin{equation}
\sum_{j=1}^{N_h} A_{ij} \alpha_j = f_i + h_i \text{ for } i=1, \dots ,N_h,
\label{eq:source}
\end{equation}
where

\begin{equation}
A_{ij} = \int_{\Omega} k(x) \text{\textbf{grad}} \phi_j \text{\textbf{grad}} \phi_i \textbf{d}x \text{ for } i, j=1, \dots ,N_h,
\label{eq:a_ij}
\end{equation}

\begin{equation}
f_i = \int_{\Omega} f \phi_i \textbf{d}x \text{ for } i=1, \dots ,N_h,
\label{eq:fi}
\end{equation}

\begin{equation}
h_i = \int_{\Gamma_N} h \phi_i \textbf{d} \Gamma \text{ for } i=1, \dots ,N_h.
\label{eq:hi}
\end{equation}

Here $A_{ij}$ are the coefficients of the stiffness matrix, and
$f_i$ and $g_i$ are the coefficients of the right-hand side
contributions from the source and the boundary terms, respectively.
All the coefficients are given by the values of the integrals in $\eqref{eq:source}$,
that are numerically approximated by a suitable quadrature rule  \cite{hughes}.

 In
order to describe this rule, let us introduce
$\widehat{\mathcal{K}}_h := \{\widehat{K}_k\}_{k=1}^{N_e}$, a
partition of the parametric domain $\widehat{\Omega}$ into $N_e$
non-overlapping subregions, that henceforth we will refer to as
$elements$. The assumptions on the parameterization $\mathbf{F}$
ensure that the physical domain $\Omega$ can be partitioned as
$$
\overline{\Omega} = \bigcup_{k=1}^{N_e}\mathbf{F}(\widehat{K}_k),
$$
and the corresponding elements $K_k:=\mathbf{F}(\widehat{K}_k)$ are
also non-overlapping. We will denote this partition by
$\mathcal{K}_h := \{K_k\}_{k=1}^{N_e}$.

For the sake of generality, let us assume that a quadrature rule is
defined on every element $\widehat{\mathcal{K}}_k$. Each of these
quadrature rules is determined by a set of $n_k \textit{ quadrature
nodes}$:
$$
\{\widehat{x}_{l,k}\} \subset \widehat{K}_k, l=1,...,n_k
$$
and by their corresponding \textit{weights}
$$
\{w_{l,k}\} \subset \mathbb{R}, l=1,...,n_k.
$$
In this work we use Gaussian quadrature rules as described in ~\cite{hughes}.

After introducing a change of variable, the integral of a generic
function $\phi \in \mathcal{L}^{1}(K_k)$ is approximated as follows
\begin{equation}
\int_{K_k} \phi dx = \int_{\widehat{K}_k} \phi (\mathbf{F}(\widehat{x}))|\det
(D\mathbf{F}(\widehat{x}))| d\widehat{x} \simeq \sum_{i=1}^{n_k} w_{l,k} \phi(x_{l,k})
|\det (D\mathbf{F}(\widehat{x}_{l,k}))|,
\end{equation}
where $x_{l,k}:=\mathbf{F}(\widehat{x}_{l,k})$ are the images of the
quadrature nodes in the physical domain, and $D\mathbf{F}$
is the Jacobian matrix of the parameterization $\mathbf{F}$.

Using the quadrature rule, the coefficients of the stiffness matrix are
numerically computed as
\begin{equation}
A_{ij} \simeq \sum_{k=1}^{N_e} \sum_{l=1}^{n_k} k(x_{l,k}) w_{l,k}
\text{\textbf{grad}} \phi_j(x_{l,k}) \text{\textbf{grad}}
\phi_i(x_{l,k}) |\det (D\mathbf{F}(\widehat{x}_{l,k}))|,
\end{equation}
while the coefficients $f_i$ of the right-hand side are approximated
as

\begin{equation}
f_i \simeq \sum_{k=1}^{N_e} \sum_{l=1}^{n_k} f(x_{l,k}) w_{l,k}
\phi_i(x_{l,k}) |\det (D\mathbf{F}(\widehat{x}_{l,k}))|.
\end{equation}

\section{Dirichlet boundary conditions}
\label{sec:dir_bc}
The imposition of homogeneous boundary conditions in IGA is
straightforward. In fact, homogeneous Neumann conditions are
satisfied naturally because they are included in the weak
formulation itself. On the contrary, as far as non-homogeneous
Dirichlet boundary conditions are concerned, a suitable
transformation has to be applied to reduce the problem to the case
of homogeneous boundary conditions since there is no place in the Galerkin formulation of the problem where the Dirichlet boundary conditions can be imposed. For this reason we must build them directly into the solution space.

In FEM the imposition of Dirichlet boundary conditions can be
performed easily since the basis functions
are interpolatory and local. On the contrary, in IGA the treatment of the Dirichlet
boundary conditions is a much more involved problem. The basis
functions (B-splines or NURBS) are not local, nor interpolatory.

In section~\ref{sec:General_IGA} we assumed that there exists a function $g_h \in \mathcal{S}_h$ such that $g_h|_{\Gamma_D} = g$, and refered to this function as a lifting. In practice, this will frequently be the case, but there will also be instances in which a lifting is only an approximation of $g$.

The standard approach consists
in constructing a function $R_g$ as suitable extension of $g$
such that
$$
R_g \in \mathcal{H}^{1}(\Omega),
$$
$$
R_g = g \text{  on  } \Gamma_D,
$$
and consider the problem $\eqref{eq:gen_weak}, \eqref{eq:var_form}$ for the difference $u - R_{g_D}$,
 see ~\cite{ciarlet_old} for details.
 To construct an approximation of
$R_{g_D}$ it is necessary to have
an approximation of the boundary data function $g$. A possible
strategy consists in obtaining an approximation $g_h$ of $g$
in the space
$$
\mathcal{S}_{h} \backslash \mathcal{V}_{h} := \langle \phi_{N_h+1},\dots,\phi_{N_h +
N_h^b} \rangle.
$$
as combination of the $N_h^b$-basis functions of $\mathcal{S}_h$ which do not
vanish on $\Gamma_D$, that is
\begin{equation}
g_h (\textbf{x}):=\sum_{N_h+1}^{N_h + N_h^b} q_i \phi_i
|_{\Gamma_D}(\textbf{x}), \text{   } \forall \textbf{x} \in
\Gamma_D.
\end{equation}
Then, an approximation of $R_g$ in $\mathcal{S}_h$, can be constructed as:
\begin{equation}
R_{g_h} (\textbf{x}):=\sum_{N_h+1}^{N_h + N_h^b} q_i \phi_i
(\textbf{x}), \text{   } \forall \textbf{x} \in \Omega_h.
\label{eq:bound1}
\end{equation}
The control variables $q_i$ in $\eqref{eq:bound1}$ have to be computed according to the approximation method used . 

The easiest
way is to apply these conditions directly to the control variables,
as
$$
g_h (\textbf{x}):=\sum_{N_h+1}^{N_h + N_h^b} g
(\overline{\mathbf{\xi}}_i) \phi_i |_{\Gamma_D}(\textbf{x}), \text{
} \forall \textbf{x} \in \Gamma_D;
$$
for proper values $\overline{\mathbf{\xi}}_i, \text{  }i=N_h+1,\dots,
N_h+N_h^b$. Such a strategy may result in a lack of accuracy in case
of inhomogeneous conditions and a poor accuracy of the
approximated solution, even if spaces $\mathcal{V}_h$ with high approximation
power are used.

A standard approach consists in imposing Dirichlet
boundary conditions by interpolation of the boundary data $g$. Another way of imposing Dirichlet boundary conditions is to use least-squares approximation of the given function $g$ by a B-spline.
We will review both methods.

\subsection{Least-squares approximation of the Dirichlet boundary conditions}
The least-squares approximation of the function $g$ consists in
solving an optimization problem. More exactly, we need to find such
coefficients $\{\widetilde{q}_i\}_{i=N_h+1}^{N_h + N_h^b}$ which
would minimize the following integral:
\begin{equation}
\min_{\{q_i\}_{i=N_h+1}^{N_h + N_h^b}} \int_{\Gamma_D}
(g(\textbf{x}) - \sum_{N_h+1}^{N_h + N_h^b} q_i \phi_i
(\textbf{x}))^2 \textbf{d} \Gamma.
\label{eq:ls}
\end{equation}

We perform the following transformations of the integral in $\eqref{eq:ls}$:
\begin{equation}
\begin{split}
&\int_{\Gamma_D}
\left( g(\textbf{x}) - \sum_{i=N_h+1}^{N_h + N_h^b} q_i \phi_i
(\textbf{x}) \right) ^2 \textbf{d} \Gamma = \\
&\int_{\Gamma_D}
(g(\textbf{x}))^2 \textbf{d} \Gamma - 2 \int g(\textbf{x}) \sum_{N_h+1}^{N_h + N_h^b} q_i \phi_i
(\textbf{x}) \textbf{d} \Gamma + \int \sum_{i,j=N_h+1}^{N_h + N_h^b} q_i q_j \phi_i (\textbf{x}) \phi_j
(\textbf{x}) \textbf{d} \Gamma.
\end{split}
\label{eq:ls1}
\end{equation}

Let us rewrite this equation in the following form:
\begin{equation}
\int_{\Gamma_D}
\left( g(\textbf{x}) - \sum_{i=N_h+1}^{N_h + N_h^b} q_i \phi_i
(\textbf{x}) \right) ^2 \textbf{d} \Gamma =
\int_{\Gamma_D}
(g(\textbf{x}))^2 \textbf{d} \Gamma - 2 \textbf{q}^T \textbf{g} +
\textbf{q}^T \textbf{M} \textbf{q},
\label{eq:ls2}
\end{equation}

where
\begin{equation}
\textbf{q}= (q_{N_h+1} \dots q_{N_h + N_h^b}),
\end{equation}

\begin{equation}
\textbf{g}= \begin{pmatrix} \int g(\textbf{x}) \phi_{N_h+1}
(\textbf{x}) \textbf{d} \Gamma \\ \vdots \\ \int g(\textbf{x}) \phi_{N_h + N_h^b}
(\textbf{x}) \textbf{d} \Gamma
\end{pmatrix},
\end{equation}
and $\textbf{M}=\{\int \phi_i (\textbf{x}) \phi_j
(\textbf{x}) \textbf{d} \Gamma\}_{i,j=N_h+1}^{N_h + N_h^b}$ is the mass matrix.

In order to minimize our target function we just differentiate it once with respect to the coefficients $q_i$ to get the closed form solution:
\begin{equation}
\textbf{q} = \textbf{M} \backslash \textbf{g}.
\end{equation}

The matrix $\textbf{M}$ being invertible since it is the Gram matrix of the basis functions of the boundary $\partial \Omega$ which is a union of NURBS or B-spline patches in the standard $\mathcal{L}^2$ inner product.

We consider the boundary of the physical domain $\Gamma$ which is an image of the boundary $\widehat{\Gamma}$ under the isogeometric mapping $\mathbf{F}$ and is itself a union of NURBS or B-spline patches with dimensions one less than the dimension of the original patch $\Omega$. So we may decompose the boundary $\Gamma = \cup_i \gamma_i$ into a union of NURBS patches and consider the minimization problem of approximating the boundary function on each of them.

After we have computed the Dirichlet degrees of freedom $\tilde{q}_i$, we can construct the desired function
\begin{equation}
R_{g_h} (\textbf{x}):=\sum_{N_h+1}^{N_h + N_h^b} \tilde{q}_i \phi_i
(\textbf{x}), \text{   } \forall \textbf{x} \in \Omega_h.
\label{eq:bound}
\end{equation}

We now solve the homogenous problem in the variational form:
\begin{equation}
a(w_h,u_h-R_{g_h}) = L(w_h), \forall w_h \in V_h,
\end{equation}
where we have only the first $N_h$ degrees of freedom to be found from this linear system.

\subsection{Quasi-interpolation of the Dirichlet boundary conditions}
Here we consider the point-wise approximation: we are looking for a B-spline function (on the given part of the boundary) which coincides with $g$ at a number of points in the spirit of collocation.

Consider a two-dimensional domain. Assume that we deal with a part of boundary $\gamma$ that corresponds to one side of the parametric domain square and there are $n$ degrees of freedom corresponding to this part of boundary. We want to take some $n$ points on the curve $\gamma$ and make the B-spline function we are looking for coincide with $g$ at these point. This gives as a square system of linear equations:

\begin{equation}
\sum_{j=1}^{n} q_j \phi_j (x_i,y_i) = g(x_i,y_i), i=1\dots n.
\end{equation}
In order to solve this system we must require its matrix to be non-singular, and this obviously depends on the choice of the points ${(x_i,y_i)}_{i=1}^n$. In our work we chose the $n$ points to be the images of the centers of the $n$ knot-spans of the parametric domain. It may be easily shown that the matrix we obtain in such a case is nonsingular. Indeed, we know that the monomials form a totally positive basis on the interval $[0,1]$ and B-splines are obtained from them by the transformations saving the totally positivity properties ~\cite{manni}, consequently the B-spline basis we have is also totally positive and the matrix we get is nonsingular.

Thus, the system now looks like:
\begin{equation}
\sum_{j=1}^{n} q_j \widehat{\phi_j} (\xi_i,\eta_i) = g(\mathbf{F}(\xi_i,\eta_i)), i=1\dots n,
\end{equation}
for the specific $\{(\xi_i,\eta_i)\}_{i=1}^n$ we choose.

One of the questions we are interested in is how to choose these points optimally for the B-spline and NURBS geometry.

Applying the same algorithm to all the sides of the domain where the Dirichlet boundary conditions are imposed we may compute all the Dirichlet degrees of freedom.

There are also other approaches if imposing the Dirichlet boundary conditions such as local least-squares, see ~\cite{gov} or Nitsche's method \cite{harari1} .

\section{Neumann boundary conditions}
Neumann boundary conditions of the form $\eqref{eq:neum_cond}$ are frequently referred to as \textit{natural boundary conditions}. This is because of the way they automatically arise in the variational statement of a problem. Let us assume for the moment without loss of generality that $k(x)=1$ and that the Dirichlet conditions have already been imposed. Multiplying by a test function and integrating leads us to 
\begin{equation}
\begin{split}
&0=\int_{\Omega} w (\triangle u + f) d \Omega \\
&= - \int_{\Omega} \nabla w \nabla u d \Omega + \int_{\Gamma} w \nabla u \cdot \mathbf{n} d \Gamma + \int_{\Omega} w f d \Omega \\
&= - \int_{\Omega} \nabla w \nabla u d \Omega + \int_{\Gamma_N} w \nabla u \cdot \mathbf{n} d \Gamma + \int_{\Omega} w f d \Omega,
\end{split}
\end{equation}
where in the third line we have used the fact that the weighting space is defined such that $w|_{\Gamma_D}=0$.

The integration by parts has naturally introduced a boundary integral over $\Gamma_N$ that refers explicitly to the condition that we would like to impose. Using the condition $\eqref{eq:neum_cond}$ we simply replace $\nabla u \cdot \mathbf{n}$ with the value we are imposing, $h$, resulting in
\begin{equation}
- \int_{\Omega} \nabla w \nabla u d \Omega + \int_{\Gamma} w h d \Gamma + \int_{\Omega} w f d \Omega = 0.
\end{equation}

Now that we have formulated a discretized analog of the original problem in the next chapter we are going to consider a variant of Schwarz Domain Decomposition Method to solve it. 

\section{Schwarz Additive Domain Decomposition Method}

\subsection{Notations}

\begin{figure}

\begin{center}
\includegraphics[width=0.5\textwidth]{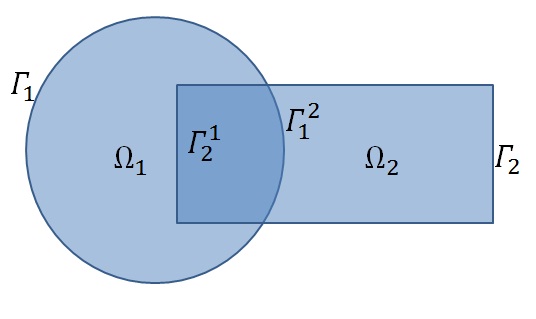}
\caption{Domain Decomposition}
\label{fig:logo_s}
\end{center}
\end{figure}

Consider the two-dimensional domain $\Omega$ as shown in the Figure \ref{fig:logo_s},
which can be represented as a union of two overlapping subdomains
$\Omega = \Omega_1 \cup \Omega_2$. In this work we consider only
domains which satisfy some reasonable conditions: they are connected
and piecewise smooth boundary. We want to solve the Poisson problem $\eqref{eq:poisson}-\eqref{eq:dir_cond}$ on this complex domain. We denote by $\Gamma_i,i=1,2$ the part of the boundary of $\Omega_i$ which belongs to the boundary $\partial \Omega$ of the whole domain and by $\Gamma_i^j,i,j=1,2,i\neq j$ the
 boundary of $\Omega_i$, which belongs to the second subdomain: $\Gamma_i^j = \partial \Omega_i \cap \Omega_j$. The earliest known domain decomposition method
is the Alternating Schwarz Method (ASM) dating back to 1870 \cite{Sc70}. The original alternating
Schwarz method is a sequential method that successively solves the
following problems: given the Dirichlet boundary conditions function $g$ on $\Gamma$ and starting with the initial guess $u_2^{0} = 0$ (in particular $u_2^{0}|_{\Gamma_1^2}=0$), iteratively obtain solutions $u_1^{n}(x,y)$ and
$u_2^{n}(x,y)$ on $\Omega_1$ and $\Omega_2$ respectively $(n \in
\mathbb{N})$ \cite{lions},

\begin{eqnarray}
   	 \triangle u_1^n = f & \quad \text{ in } \Omega_1, \label{eq:SADD1-1}\\
	u_1^n = g_1^n & \quad \text{ on }\Gamma_{1D},
\end{eqnarray} 
where $g_1^n|_{\Gamma_1^2} = u_2^{n-1}|_{\Gamma_1^2}, g_1^n|_{\Gamma_1} = g|_{\Gamma_1}$ 
and

\begin{eqnarray}
   	 \triangle u_2^n = f & \quad \text{ in } \Omega_2,\\
	u_2^n = g_2^n & \quad \text{ on }\Gamma_{2D}, \label{eq:SADD2-3}
\end{eqnarray}
where $g_2^n|_{\Gamma_2^1} = u_1^{n}|_{\Gamma_2^1}, g_2^n|_{\Gamma_2} = g|_{\Gamma_2}$.
It is important to notice that the solution of the first problem is
required before the second problem can be solved. The part of the boundary
conditions of the second problem, $g_2^n = u_1^{n}|_{\Gamma_2^1}$, 
is the current solution of the first problem on
the interface. It is not
difficult to see that in this simple case of two subdomains it is an analog of a \textit{block Gauss-Seidel} method and the two
problems cannot be solved in parallel. 
For a straight forward
parallelization of the algorithm an analog of the \textit{block Jacobi} method is applied giving us the Schwarz Additive method:

\begin{eqnarray}
   	 \triangle u_1^n = f & \quad \text{ in } \Omega_1,\label{eq:SADD3-1}\\
	u_1^n = g_1^n & \quad \text{ on }\Gamma_{1D},
\end{eqnarray}
where $g_1^n|_{\Gamma_1^2} = u_2^{n-1}|_{\Gamma_1^2}, g_1^n|_{\Gamma_1} = g|_{\Gamma_1}$, and
\begin{eqnarray}
   	 \triangle u_2^n = f & \quad \text{ in } \Omega_2,\\
	u_2^n = g_2^n & \quad \text{ on }\Gamma_{2D}, \label{eq:SADD4-3}
\end{eqnarray}
where $g_2^n|_{\Gamma_2^1} = u_1^{n-1}|_{\Gamma_2^1}, g_2^n|_{\Gamma_2} = g|_{\Gamma_2}$.

The iterative procedure starts given two initial guesses
$u_1^{0}$ and $u_2^{0}$ on each subdomain. 
The boundary condition of one of the subproblems,
 $g_j^n = u_i^{n-1}|_{\Gamma_j^i}, i,j=1,2, i \neq j$ is now the solution of the
previous iterate on the interface, hence the subproblems can be
solved independently. In the general case of many
subdomains a great deal of parallelism can be introduced, using a
method called multi-coloring. This is a graph theory technique that
identifies each subdomain with a color such that disjoint subdomains
have the same color; usually only a small number of colors is
needed. Thus, the subproblems defined on subdomains of the same
color can be solved independently in parallel and only the number of
distinct colors different threads is needed to perform the computations. We are using such techniques below when considering 3-dimensional multi-patched domains.

For the ASDDM as presented here, one can see
that there is no notion of a global approximate solution. One can
consider defining the global solution by choosing the subdomain iterative solution on each
subdomain and any weighted average within the overlap:
$$
u^{n} = \chi_1 v^{n} + \chi_2 w^{n},
$$
where $\chi_1 = 1$ on $\Omega_1 \backslash (\Omega_1 \cap
\Omega_2)$,  $\chi_2 = 1$ on $\Omega_2 \backslash (\Omega_1 \cap
\Omega_2)$ and $\chi_1 + \chi_2 = 1$ on $\Omega_1 \cap \Omega_2$. This approach extends easily to a case of multi-patched domain decomposition.

\subsection{Convergence conditions}

The iterations are performed until certain convergence conditions are met. Throughout this work we consider the $\mathcal{L}^2$ errors as the termination condition. In the model problems when the exact solution is known the errors are calculated as a $\mathcal{L}^2$ norms between the approximate solution and the exact one. 

Assume the exact solution $u_{ex} \in \mathcal{H}^1$ is given. Then the error is calculated as
\begin{equation}
\begin{split}
&||u_h - u_{ex}||_2^2 = \int_{\Omega}(u_h - u_{ex})^2 dx = \int_{\Omega}\left( \sum_{j \in \mathcal{J}_i \cup \mathcal{J}_b} \alpha_j \phi_j - u_{ex} \right)^2 dx = \\
&\sum_{k} \int_{K_k} \left( \sum_{j \in \mathcal{J}_i \cup \mathcal{J}_b} \alpha_j \phi_j - u_{ex} \right)^2 dx. 
\end{split}
\end{equation}

We sum up the integrals on the elements and perform numerical integration element-wise. Assume that an integration rule is given by the mesh of points $\{\widehat{x}_{l,k}\}_{l,k}$ and weights $\{w_{l,k}\}_{l,k}$. For each element $K_{k}$ we obtain:
\begin{equation}
\begin{split}
& \int_{K_k} \left( \sum_{j \in \mathcal{J}_i \cup \mathcal{J}_b} \alpha_j \phi_j - u_{ex} \right)^2 dx \\
& = \int_{\widehat{K}_k}  \left( \sum_{j \in \mathcal{J}_i \cup \mathcal{J}_b} \alpha_j \phi_j (\mathbf{F}(\widehat{x}))) - u_{ex}(\mathbf{F}(\widehat{x}))) \right)^2 \det
(D\mathbf{F}(\widehat{x}))| d\widehat{x} \\
&\simeq \sum_{i=1}^{n_k} w_{l,k} \left( \sum_{j \in \mathcal{J}_i \cup \mathcal{J}_b} \alpha_j \phi_j (\mathbf{F}(\widehat{x}_{l,k}))) - u_{ex}(\mathbf{F}(\widehat{x}_{l,k}))) \right)^2 
|\det (D\mathbf{F}(\widehat{x}_{l,k}))|.
\end{split}
\end{equation}

This numerical integration formula was implemented in the code we used. An analogous formula for the $\mathcal{H}^1$ errors can also be developed.

\section{Weak form of ASDDM}

We may now obtain a variational form of the ASDDM method for the equations $\eqref{eq:SADD3-1}-\eqref{eq:SADD4-3}$. 

We are looking for two functions $u_i \in \mathcal{H}^1(\Omega_i), i=1,2$ such that their weighted sum $u =
\chi_1 u_1 + \chi_2 u_2$ satisfies the Poisson equation $\eqref{eq:poisson}-\eqref{eq:dir_cond}$ or its weak form $\ref{eq:var_form}$.

This  weak form of the iterative ASDDM now reads as:

\begin{algorithm}
\caption{ ASDDM Weak form}
\begin{algorithmic}
\STATE  Given $u_i^0 \in \mathcal{H}^1(\Omega_i), i=1,2$, such that $u_i^0|_{\Gamma_i} =  g|_{\Gamma_i}, i=1,2 $, 
\STATE Define convergence  level $\eps$
\WHILE{ $Error \ge  \eps $ }
\STATE \emph{Find $u_i^{n} \in \mathcal{H}^1(\Omega_i):$ such that
$u_i^n|_{\Gamma_i^j} = u_j^{n-1}|_{\Gamma_i^j} \text{ and }
u_i^n|_{\Gamma_i} = g|_{\Gamma_i}$}
\STATE
\emph{$$ \int_{\Omega_i} \textbf{grad} v \cdot \textbf{grad} u_i^n
\textbf{d} \Omega = \int_{\Omega_i} v f \textbf{d} \Omega +
\int_{\Gamma_{iN}} v h \textbf{d} \Gamma $$}
 \STATE \emph{ for any $v \in \mathcal{H}_0^1(\Omega_i),i,j=1,2, i \neq j$. }
 \ENDWHILE
 \end{algorithmic}
\end{algorithm}

\section{Operator form of the algorithm}
In order to proceed to the finite dimensional space computations from this continuous settings we need to introduce a few operators. Assume that the solution $u_i^n$ on the $\Omega_i$ domain is sought in the finite dimensional space $\mathcal{S}_i$. Given a function $u_j^n \in \mathcal{H}^1(\Omega_j)$ our aim is to impose the Dirichlet boundary conditions on the $\Omega_i$ domain, which requires us to find a function $g_i \in \mathcal{S}_i$ which will be close in some specific sense to the function $u_j$ on the interface boundary $\Gamma_i^j$.

We proceed as following: first we need to project the function $u_j^n$ onto the boundary $\Gamma_i^j$ to get a function $\widetilde{g} \in \mathcal{L}^2(\Gamma_i^j)$ representing the Dirichlet boundary conditions on this part of the boundary. We refer to this operator as a \textit{trace operator}.

The next steps depend on whether the space $\mathcal{S}_i$ contains a function which is an extension of $\widetilde{g}$ or not. If it does we just take any extension of $\widetilde{g}$ satisfying the Dirichlet boundary conditions on the parts of the boundary $\Gamma_i$ and continue the solution process as in the regular IGA case, see section (3.2.1). We refer to the operator building the extension as \textit{extension operator}. If $\mathcal{S}_i$ does not contain an extension of $\widetilde{g}$ we need to find a function in $\mathcal{S}_i$ that would approximate the $\widetilde{g}$ on the boundary $\Gamma_i^j$. We perform this approximation in two steps. First, we build an auxilary space $\mathcal{S}_{\Gamma_i^j}$ of functions that are just restrictions of the functions in $\mathcal{S}_i$ onto the boundary $\Gamma_i^j$. We then find among these functions the one $\widetilde{g}_i \in {\mathcal{S}_i}|_{\Gamma_i^j}$ that is closest to the $\widetilde{g}$. Operator doing this is called \textit{approximation operator}. The process is then continued as before by applying the extension operator to $\widetilde{g}_i$ to obtain $g_i \in \mathcal{S}_i$.

We now define the operators we have just introduced and rewrite the iterative algorithm in the operator notations.

Loosely speaking the trace operator is defined as following: given a function $w$ from a smooth enough functional space defined on $\Omega_j$ the trace operator returns function defined on $\Gamma_i^j$ which coincides with $w$ on  $\Gamma_i^j, i \neq j$.

In order to define the trace operator exactly we use the following
\begin{thm}(Trace theorem)
Let $\Omega \subset \mathbb{R}^n$ be bounded with piecewise smooth boundary. Then there exists a bounded linear mapping
\begin{equation}
\Pi: \mathcal{H}^1(\Omega) \rightarrow \mathcal{L}^2(\Omega),
\end{equation}
\begin{equation}
||\Pi (u)||_{0,\partial \Omega} \leq C ||u||_{1,\Omega},
\end{equation}
such that $\Pi (u) = u|_{\partial \Omega} \text{ } \forall u \in C^1(\overline{\Omega})$.
\end{thm}
The proof can be found in ~\cite{lm}.

Since we assume the domain boundaries to be piecewise smooth and the solutions belong to $\mathcal{H}^1(\Omega_j)$ spaces, the Trace theorem provides us with the definition of the trace operator:

$$
P_{i} : \mathcal{H}^1(\Omega_j) \rightarrow \mathcal{L}^2(\Gamma_i^j)
$$
$$
u_j \rightarrow u_j|_{\Gamma_i^j}.
$$

Now assume that $\mathcal{S}_i = \mathcal{H}^1(\Omega_i)$ and the extension is possible. Given a function defined on $\Gamma_i^j$ we need to extend this function to the whole subdomain $\Omega_i$ while keeping the smoothness properties (almost everywhere with respect to the Lebesgue measure) of the function and in such a way that the new function will coincide with the Dirichlet boundary condition $g$ on the outer boundary $\Gamma_i$.
More precisely, the extension operators look like:
$$
E_i : \mathcal{H}^1(\Gamma_i^j) \rightarrow \mathcal{H}^1(\Omega_i)
$$
$$
v_i \rightarrow u_i \text{ such that } u_i|_{\Gamma_i^j} = v_i, u_i|_{\Gamma_i} = g|_{\Gamma_i}.
$$
Here we assume that the conditions $ u_i|_{\Gamma_i^j} = v_i, u_i|_{\Gamma_i} =
g|_{\Gamma_i}$ define some continuous $\mathcal{H}^1(\Omega_i)$ function. This condition is satisfied throughout the whole discussion.

Denote the corresponding billinear and linear forms on the subdomains as $a_i, L_i, i=1,2.$ The continuous version of the ASDDM now reads as follows.

\begin{algorithm}
\caption{ ASDDM  with Trace Operator}
\begin{algorithmic}
\STATE  Given $u_i^0 \in \mathcal{H}^1(\Omega_i), i=1,2$, such that $u_i^0|_{\Gamma_i} =  g|_{\Gamma_i}, i=1,2 $, 
\STATE Define convergence  level $\eps$
\WHILE{ $Error \ge  \eps $ }
\STATE \emph{Find $u_i^{n} \in \mathcal{H}^1(\Omega_i):$ such that}
\emph{$$a_i(u_i^n - E_i P_i u_j^{n-1},v_i) = L_i (v_i) - a_i(E_i P_i u_j^{n-1},v_i)$$} 
\STATE \emph{for any $v_i \in \mathcal{H}_0^1(\Omega_i), i,j=1,2, i \neq j.$}
 \ENDWHILE
 \end{algorithmic}
\end{algorithm}

When we consider any iterative procedure the main question that arises is whether this process converges, and under which conditions. The convergence of this classical Schwarz algorithm was shown by P.-L. Lions ~\cite{lions} under some regularity assumptions.

In the next section we consider the case when the approximation need to be performed before we can extend the boundary function to the whole domain and give a deeper overview of the trace operators we implemented.

\section{Discretized Schwarz Additive Algorithm and Trace Operators}

We now assume that we deal with the finite dimensional function spaces $\mathcal{S}_1$ and $\mathcal{S}_2$ instead of $\mathcal{H}^1(\Omega_1)$ and $\mathcal{H}^1(\Omega_2)$, respectively. It means that
we can not just project the function from $\mathcal{S}_j$ onto the boundary
$\Gamma_i^j$ and extend it to the whole $\Omega_i$ domain, since this projection, generally speaking, will not
belong to $\mathcal{S}_i$ restricted to $\Gamma_i^j$. So we have to
introduce an additional approximation operator as was explained in the previous chapter. This operator acts on functions from $\mathcal{H}^1(\Gamma_i^j)$ and maps them into the subspace ${\mathcal{S}_i}|_{\Gamma_i^j} \subset \mathcal{H}^1(\Gamma_i^j)$ of restrictions of functions from $\mathcal{S}_i$ onto the boundary $\Gamma_i^j$.    

More specifically, in our IGA environment we want this operator to build a B-spline or NURBS approximation of the function in hand. In the previous notations the approximation operator will look as following:
$$
A_i : \mathcal{H}^1(\Gamma_i^j) \rightarrow {\mathcal{S}_i}|_{\Gamma_i^j}
$$
$$
w \rightarrow v\in \mathcal{S}|_{\Gamma_i^j}, \text{approximating } w,
$$
where the approximating function can be thought of as minimizing distance to the function being approximated in some specific norm.

With this notations the numerical algorithm will look like:

\begin{algorithm}
\caption{ Discrete ASDDM  with Trace Operator}
\begin{algorithmic}
\STATE  Given initial guesses $u_i^0 \in \mathcal{S}_i,i=1,2$ such that
$u_i^0|_{\Gamma_i} = g_h|_{\Gamma_i},i=1,2$
\STATE Define convergence  level $\eps$
\WHILE{ $Error \ge  \eps $ }
\STATE \emph{Find $u_i^{n} \in \mathcal{S}_i$ such that :
$$a_i(u_i^n - E_i A_i P_i u_j^{n-1},v_i) = L_i(v_i) - a_i(E_i A_i P_i u_j^{n-1},v_i)$$}
\STATE \emph{for any $v_i \in \mathcal{V}_i, i,j=1,2, i \neq j.$}
 \ENDWHILE
 \end{algorithmic}
\end{algorithm}

%
%

\begin{figure}
\begin{center}
\includegraphics[width=0.9\textwidth]{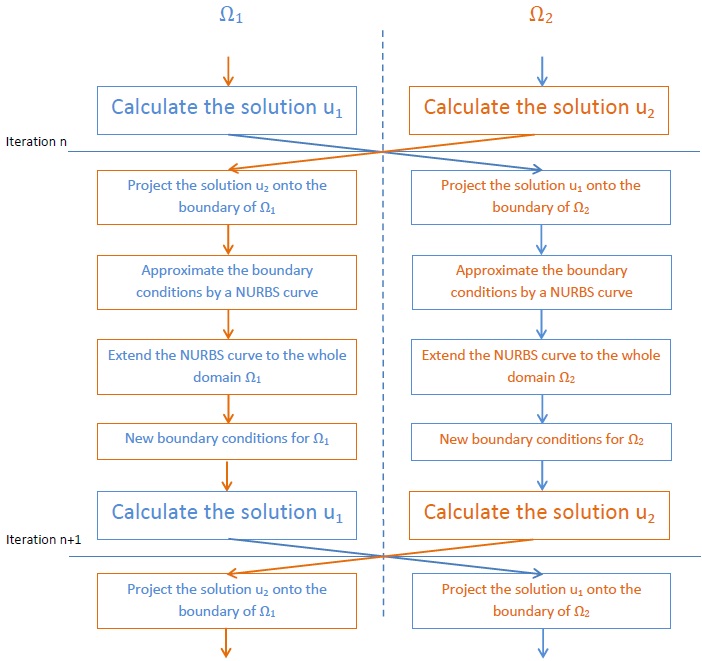}
\caption{Operations performed on each iterative step in ASDDM with two subdomains $\Omega_1$ and $\Omega_2$}
\label{fig:algo_struct}
\end{center}
\end{figure}

Figure \ref{fig:algo_struct} illustrates the described ASDDM applied to a two domains problem.

\section{The Trace Operator}

In the present work we consider two ways of projection and compare their performance.
Assume that both domains $\Omega_1$ and $\Omega_2$ are B-spline surfaces with the corresponding piecewise invertible mappings $\mathbf{F}_1:\widehat{\Omega}_1 \rightarrow \Omega_1$ and $\mathbf{F}_2:\widehat{\Omega}_2 \rightarrow \Omega_2$, where the parametric domains are $\widehat{\Omega}_1  = \widehat{\Omega}_2 = [0,1] \times [0,1]$. 

\subsection{Exact trace operator}
We define the exact projetion operator as before:
$$
P^e_i : \mathcal{H}^1(\Omega_j) \rightarrow \mathcal{H}^1(\Gamma_i^j)
$$
$$
u_j \rightarrow u_j|_{\Gamma_i^j}.
$$

After each iteration on each domain the approximate solution we get
is a linear combination of images of B-spline (or NURBS) basis functions. Let
us consider the solution $u_j^{n-1}$ obtained after $n-1$-th
iteration on the subdomain $\Omega_j$:
\begin{equation}
u_j^{n-1}(x, y)=\sum_{k \in \mathcal{J}} \phi_k (x, y)
{u_j}_k^{n-1},
\end{equation}
where  $\mathbf{u_j} = \{{u_j}_i\}_{i \in \mathcal{J}}$ is the
vector of degrees of freedom.

We want to evaluate this function at some point $(x,y) \in \Omega_j$. Recall that the pull-back $\mathbf{F}_j$ always exists since the geometrical mapping $\mathbf{F}_j$ that we use was chosen to be invertible: $(\xi,\eta) =
\mathbf{F}_j^{-1}(x,y)$. The problem here is that we cannot obtain a closed analytic form of $\mathbf{F}_j^{-1}$, which means that in the real computations we will need to use some numerical methods to get the pre-image $(\xi,\eta)$ of the point $(x,y)$.

\begin{figure}
\begin{center}
\includegraphics[width=0.8\textwidth]{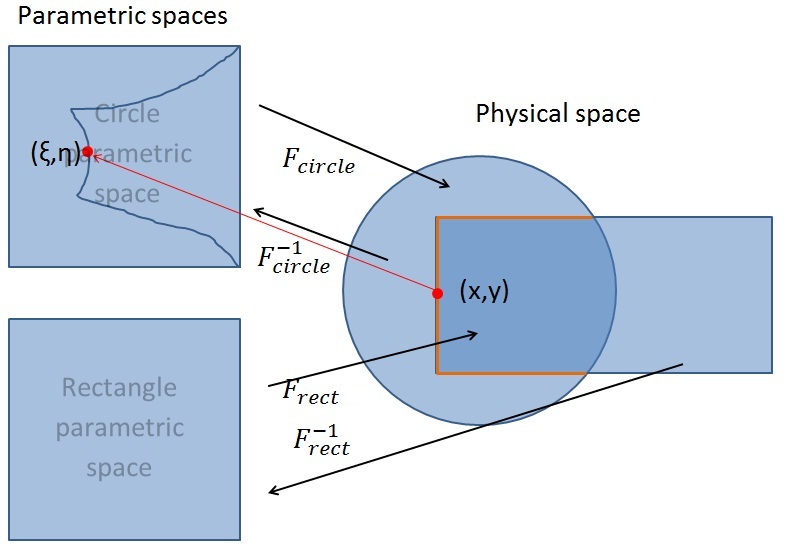}
\label{fig:pullback2}
\caption{Domains and the pullback $F_{circle}^{-1}(\Gamma_2^1)$ of the boundary $\Gamma_2^1$}
\end{center}
\end{figure}

The resulting function $P^e_i(u_j^{n-1})$ is obtained by evaluating the solution $u_j^{n-1}$ at any point $(x,y)$ in this way.

In the Figure ~\ref{fig:pullback2} you can see an example of a pullback of the boundary $\Gamma_2$ into the parametric space $\widehat{\Omega}_1$ 

\subsection{Interpolation trace operator}

The second approach of performing the projection we used was the interpolation trace operator. Given the geometrical mapping $\mathbf{F}$ we can generate a mesh $\mathcal{M}$ in the physical space by
applying this mapping to some mesh $\widehat{\mathcal{M}}$ in the parametric space. On the $n$-th iteration we can evaluate the solution $u_j^n,  j=1,2$,  in
the physical domain $\Omega_j$ at the mesh points we have just constructed.
Our aim is to approximate this solution at any point of the domain $\Omega_j$. We obtain it by taking the values of the function $u_j^n$ in the constructed mesh points and interpolating them into the points of interest, e.g. points on the curve $\Gamma_i^j$. 

\begin{figure}
\begin{center}
\includegraphics[width=0.4\textwidth]{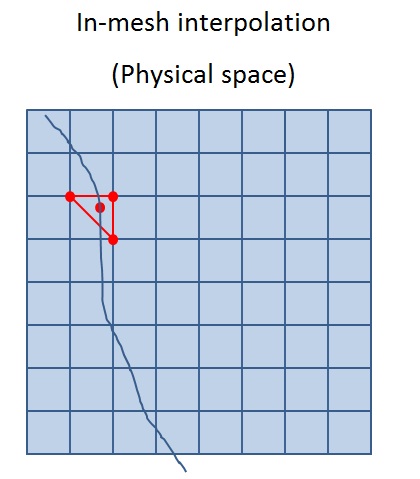}
\label{fig:int_proj}
\caption{Interpolation projection}
\end{center}
\end{figure}

There is a wide choice of interpolation methods that we can use. In
our research we used linear and cubic interpolations. The
interpolation is obtained in the following way: we, first, divide the physical domain
into triangles by applying the Delaunay triangulation algorithm to the mesh $\mathcal{M}$. Then for any point $(x,y) \in \Gamma_i^j$ we define to which triangle it belongs and 
interpolate the value of the function $u_j^n$ at this point using
the values of this function in the vertices of the triangle we have chosen.
$$
P_i^j : \mathcal{H}^1(\Omega_j) \rightarrow \mathcal{H}^0(\Gamma_i^j)
$$
$$
u_j \rightarrow I_{\mathcal{D}}(u_j),
$$
where $I_{\mathcal{D}}$ is the interpolation operator, depending on the triangulation $\mathcal{D}$ of the mesh $\mathcal{M}$.
We will use the notation $P^l$ to denote the linear interpolation trace operator.

In the next section we treat a one-dimensional example of ASDDM. We give an overview on the construction of the solver matrices and their structure. 

\section{One-dimensional ASDDM : Discretized Algorithm }

\begin{figure}
\begin{center}
\includegraphics[width=0.85\textwidth]{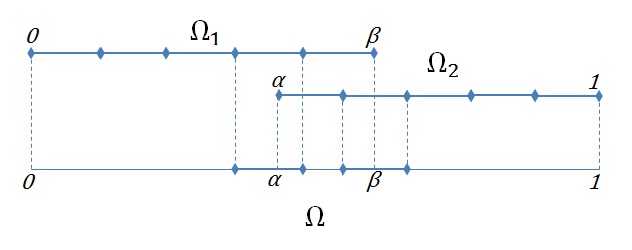}
\label{fig:1d_dec}
\caption{One-dimentional Domain Decomposition}
\end{center}
\end{figure}

To illustrate the ASDDM on  non-matching grids we consider here a simple one-dimensional Poisson equation on a line segment domain $\Omega = [0,1]$. 
\begin{equation}
-u_{xx} = f \text{ on } \Omega_1, u(0) = u(1) = 0.
\end{equation}

We decompose the domain into two overlapping domains $\Omega_1 = [0,\beta) \text{ and } 
\Omega_2 =
(\alpha,1]$,  $\alpha < \beta$ such that $\Omega_1 \cup \Omega_2 = \Omega = [0,1]$ as shown in the Figure ~\ref{fig:1d_dec}. We assume for convenience that the mapping F
is the identity mapping, which implies that we
may consider all the equations in the parametric spaces.
We denote the knot vectors as $\Xi_1 = [0, \dots \beta]$ and $\Xi_2 =
[\alpha, \dots 1]$ and the numbers of degrees of freedom as $N_{h1}$ and $N_{h2}$, respectively. 
Also we assume that the knot
vectors are uniform and open and denote the degrees of the
polynomials $p_1$ and $p_2$ on the subdomains $\Omega_1$ and
$\Omega_2$, respectively. To avoid too many indices we denote the solution on $\Omega_1$ as $v$
and the solution on $\Omega_2$ as $w$.

The two equations which we are going to
discretize are:
\begin{equation}
v^{n}_{xx} = f \text{ on } \Omega_1, v^{n}(0) = 0, v^{n}(\beta) =
w^{n-1}(\beta).
\end{equation}
\begin{equation}
w^{n}_{xx} = f \text{ on } \Omega_2, w^{n}(0) = 0, w^{n}(\alpha) =
v^{n-1}(\alpha).
\end{equation}

The numbering of the basis functions is exactly the same as we
used in section ~\ref{sec:General_IGA} .

Consider now the first subdomain $\Omega_1$. Only the first basis function $\phi_1$
is non-zero at $0$ satisfying $\phi_1(0)=1$. In the similar way at the second end of $\Omega_1$ $\phi_{N_{h1}}(\beta) = 1$ and all the other basis functions vanish at point $\beta$. It means, that in the one-dimensional
case the imposition of Dirichlet boundary conditions becomes
trivial: we just need to take the first and the last degrees of
freedom equal to the projections of the solution on the other
subdomain at the boundary points; the approximation operators $A_i$ become identities. 
\\
Denote the approximate solution on the first subdomain as $\widetilde{v}$ and its vector representation in the B-spline basis as $\mathbf{v} = \{ v_i \}_{i=1}^{N_{h1}}$. Since the boundary
condition at $0$ is constantly zero, the first degree of freedom is $0$: $v_1 = 0$. The treatment of the boundary condition at the
point $\beta$ is more involved. We need to project the solution
$\widetilde{w}^{n-1}$ of the previous iteration on the second subdomain
$\Omega_2$ onto the boundary of $\Omega_1$ which in our case consists
of the only point $\beta$. So we can write $v_{N_{h1}}^n =
P_1(\widetilde{w}^{n-1}(\beta)),$ where $P_1$ is the trace operator.

The same reasoning applies on $\Omega_2$.

Denote the vector of degrees of
freedom of the solution $\widetilde{w}$ as $\mathbf{w}$.
The solution $\widetilde{w}^{n-1}$ now reads as
$$
\widetilde{w}^{n-1} = \sum_{i=1}^{N_{h2}} \psi_i w_i^{n-1}.
$$
All the operators are acting on finite dimensional spaces, thus they have matrix representations in the basis we have chosen. Since the degree of B-splines on $\Omega_2$
is $p_2$ we may conclude that there are no more than $p_2+1$
non-zero basis functions $\psi_i$ at the point $\beta$. As we
discussed above, we may use different trace operators in order
to project the solution on the subdomain $\Omega_2$ onto the
boundary $\Gamma_1^2$ of the subdomain $\Omega_1$. In what
follows we will consider both cases: the exact trace operator
and the interpolation trace operator.

We consider now the two kinds of trace operators we introduced in the previous chapter.

\subsection{Exact trace operator $P^e$} 

We take the value $v_{N_{h1}}^{n}$ of the boundary degree of
freedom to be equal exactly to the value of the function
$\widetilde{w}^{n-1}(\beta)$. For this trivial trace operator, the
value $v_{N_{h1}}^{n}$ we get is:
$$
v_{N_{h1}}^{n} = (0, \dots, \psi_i(\beta),
\psi_{i+1}(\beta),
\dots, \psi_{i+p_2}(\beta), \dots, 0) \cdot \begin{pmatrix} 0 \\ \vdots \\ w_{i}^{n-1} \\ w_{i+1}^{n-1} \\
\vdots \\ w_{i+p_2}^{n-1} \\ \vdots \\ 0
\end{pmatrix}
$$

The operator matrix for the boundary of the $\Omega_1$ subdomain is:
$$
P^e_1 =
\begin{pmatrix}
0 & \dots & 0 & 0 & \dots & 0 & \dots & 0 \\
0 & \dots & 0 & 0 & \dots & 0 & \dots & 0 \\
\vdots & \dots & \vdots & \vdots & \dots & \vdots & \dots & \vdots \\
0 & \dots & \psi_i(\beta) & \psi_{i+1}(\beta) & \dots &
\psi_{i+p_2}(\beta) & \dots & 0
\end{pmatrix},
$$
where $P^e_1 \in \mathbb{R}^{N_{h1} \times N_{h2}}.$

Denote by $A_1$ the stiffness matrix for the first subdomain.
As before, we partition the sets of indices of basis functions $\mathcal{J} = \{1,2,\dots,N_{hj}\},j=1,2$ into two subsets. The subset $\mathcal{I} \subset \mathcal{J}$ of the inner degrees of freedom and $\mathcal{B} \subset \mathcal{J}$ of the boundary degrees of freedom. In the one dimensional case $\mathcal{B} = \{1, N_{hj}\}$ and $\mathcal{I} = \mathcal{J} \backslash \{1, N_{hj}\} = \{2,3,\dots,N_{hj}-1\}, j=1,2.$ So
the restriction of the stiffness matrix corresponding to the inner degrees
of freedom reads as $A_1(\mathcal{I}, \mathcal{I})$. 

Since the basis functions are not local, when we impose the Dirichlet boundary
conditions on some of the degrees of freedom, we need to subtract
the function $v_{N_{h1}}^{n} \phi(\xi)$ from the solution
$\widetilde{v}^{n}$. We do it by subtracting the corresponding values from the
degrees of freedom of the basis functions
$\{\phi\}_{i=N_{h1}-k}^{N_{h1}}$ which domain intersects with the
domain of $\phi_{N_{h1}}$. Now the discretized equation for the subdomain $\Omega_1$ can be rewritten as:
\begin{equation}
\begin{split}
&\begin{pmatrix}
1 & 0 &\dots & 0 & 0 \\
0 & & A_1 (\mathcal{I},\mathcal{I}) & & 0 \\
0 & 0 & \dots & 0 & 1
\end{pmatrix} \cdot \mathbf{v}^{n} =
\begin{pmatrix} 0 \\ \mathbf{f_1} \\ 0
\end{pmatrix} +
\begin{pmatrix}
1 & \dots & 0 \\
& - A_1 (\mathcal{I},\mathcal{B})& \\
0 & \dots & 1
\end{pmatrix} \cdot\\
& \cdot
\begin{pmatrix}
0 & \dots & 0 & 0 & \dots & 0 & \dots & 0 \\
0 & \dots & 0 & 0 & \dots & 0 & \dots & 0 \\
\vdots & \dots & \vdots & \vdots & \dots & \vdots & \dots & \vdots \\
0 & \dots & \psi_i(\beta) & \psi_{i+1}(\beta) & \dots &
\psi_{i+p_2}(\beta) & \dots & 0
\end{pmatrix} \cdot \mathbf{w}^{n-1},\\
\end{split}
\label{eq:matrix_form_1}
\end{equation}
where the vector $\mathbf{f_1}$ corresponds to the inner degrees of
freedom of the first subdomain. Exactly the same reasoning can be
applied to the second subdomain $\Omega_2$ to get the following discretized
equation:
\begin{equation}
\begin{split}
&\begin{pmatrix}
1 & 0 &\dots & 0 & 0 \\
0 & & A_2 (\mathcal{I},\mathcal{I}) & & 0 \\
0 & 0 & \dots & 0 & 1
\end{pmatrix} \cdot \mathbf{w}^{n} =
\begin{pmatrix} 0 \\ \mathbf{f_2} \\ 0
\end{pmatrix} +
\begin{pmatrix}
1 & \dots & 0 \\
& - A_2 (\mathcal{I},\mathcal{B})& \\
0 & \dots & 1
\end{pmatrix} \cdot\\
& \cdot
\begin{pmatrix}
0 & \dots & \phi_j(\alpha) & \phi_{j+1}(\alpha) & \dots &
\phi_{j+p_1}(\alpha) & \dots & 0 \\
\vdots & \dots & \vdots & \vdots & \dots & \vdots & \dots & \vdots \\
0 & \dots & 0 & 0 & \dots & 0 & \dots & 0 \\
0 & \dots & 0 & 0 & \dots & 0 & \dots & 0 \\
\end{pmatrix} \cdot \mathbf{v}^{n-1}.\\
\end{split}
\label{eq:matrix_form_2}
\end{equation}

Let us denote the vector of degrees of freedom of the \lq\lq{}whole\rq\rq{} approximate solution
$\widetilde{u}^n$ as $\mathbf{u}^n = \begin{pmatrix} \mathbf{v}^{n} \\ \mathbf{w}^{n}
\end{pmatrix}$. Now, we can unite the two equations \eqref{eq:matrix_form_1} and \eqref{eq:matrix_form_2} into
a single matrix equation
\begin{equation}
A \cdot \mathbf{u}^n = f + A_{dir} \cdot P \cdot \mathbf{u}^{n-1},
\end{equation}
where
$$ A =
\begin{pmatrix}
\begin{matrix}
1 & \dots & 0 \\
0 & A_1 (\mathcal{I},\mathcal{I}) & 0\\
0 & \dots & 1\\
\end{matrix} & O \\
O &
\begin{matrix}
1 & \dots & 0\\
0 & A_2 (\mathcal{I},\mathcal{I}) & 0 \\
0 &\dots & 1 \\
\end{matrix}
\end{pmatrix}
= 
\begin{pmatrix}
\tilde{A}_1 & O \\
O & \tilde{A}_2
\end{pmatrix},
$$
with $\tilde{A}_1 \in \mathbb{R}^{N_{h1} \times N_{h1}},\tilde{A}_2 \in \mathbb{R}^{N_{h2} \times N_{h2}}, $ 
$$ f =
\begin{pmatrix}
0 \\ \mathbf{f_1} \\ 0 \\ 0 \\ \mathbf{f_2} \\0 \\
\end{pmatrix},
$$
$$ A_{dir} =
\begin{pmatrix}
\begin{matrix}
1 & \dots & 0 \\
& - A_1 (\mathcal{I},\mathcal{B}) & \\
0 & \dots & 1 & \\
\end{matrix} & O\\
O &
\begin{matrix}
1 & \dots & 0\\
& -A_2 (\mathcal{I},\mathcal{B}) & \\
0 &\dots & 1 \\
\end{matrix}
\end{pmatrix}
= 
\begin{pmatrix}
\widetilde{A}_{dir1} & O \\
O & \widetilde{A}_{dir2}
\end{pmatrix},
$$
with $\widetilde{A}_{dir1} \in \mathbb{R}^{N_{h1} \times N_{h1}},\widetilde{A}_{dir2} \in \mathbb{R}^{N_{h2} \times N_{h2}}, $

\begin{equation}
\begin{split}
& P =
\begin{pmatrix}
O &
\begin{matrix}
0 & \dots & 0 & \dots & 0 & \dots & 0 \\
0 & \dots & 0 & \dots & 0 & \dots & 0 \\
\vdots & \dots & \vdots & \dots & \vdots & \dots & \vdots \\
0 & \dots & \psi_i(\beta) & \dots &
\psi_{i+p_2}(\beta) & \dots & 0 \\
\end{matrix}\\
\begin{matrix}
0 & \dots & \phi_j(\alpha) & \dots &
\phi_{j+p_1}(\alpha) & \dots & 0 \\
\vdots & \dots & \vdots & \dots & \vdots & \dots & \vdots \\
0 & \dots & 0 & \dots & 0 & \dots & 0 \\
0 & \dots & 0 & \dots & 0 & \dots & 0 \\
\end{matrix} & O
\end{pmatrix}\\
&
= 
\begin{pmatrix}
O & \tilde{P}_1 \\
\tilde{P}_2 & O
\end{pmatrix},
\end{split}
\end{equation}
with $P_1^e \in \mathbb{R}^{N_{h1} \times N_{h2}},P_2^e \in \mathbb{R}^{N_{h2} \times N_{h1}} .$

The iterative scheme is given by:
\begin{equation}
\begin{split}
& \begin{pmatrix}
\tilde{A}_1 & O \\
O & \tilde{A}_2
\end{pmatrix}
\cdot \mathbf{u}^n = f +\begin{pmatrix}
\widetilde{A}_{dir1} & O \\
O & \widetilde{A}_{dir2}
\end{pmatrix}
\cdot 
\begin{pmatrix}
O & P_1^e \\
P_2^e & O
\end{pmatrix}
\cdot \mathbf{u}^{n-1} = \\
&f+\begin{pmatrix}
O & \widetilde{A}_{dir1} P_1^e \\
\widetilde{A}_{dir2} P_2^e & O
\end{pmatrix}
\cdot \mathbf{u}^{n-1},
\end{split}
\end{equation}

\begin{figure}
\begin{center}
\includegraphics[width=01\textwidth]{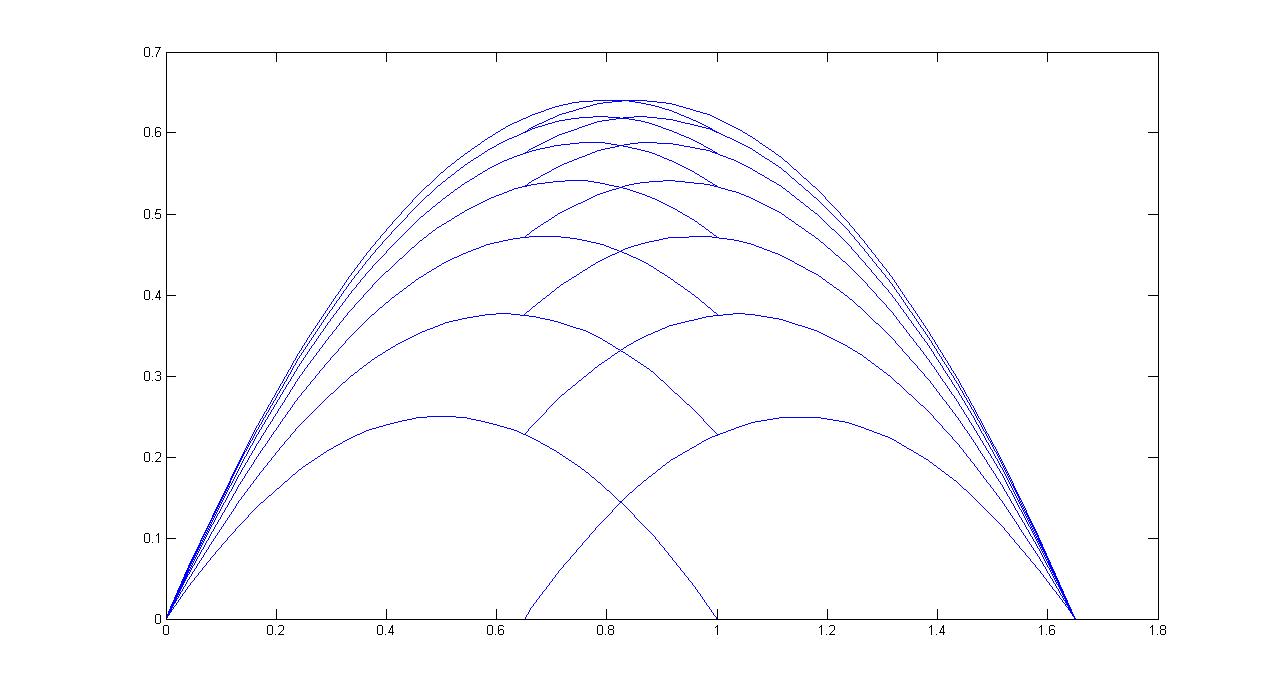}
\caption{One-dimensional Schwarz Method for the equation $u_{xx}=-1$, using the exact trace operator $P^e$ and zero initial guesses}
\end{center}
\end{figure}

\begin{figure}
\begin{center}
\includegraphics[width=01\textwidth]{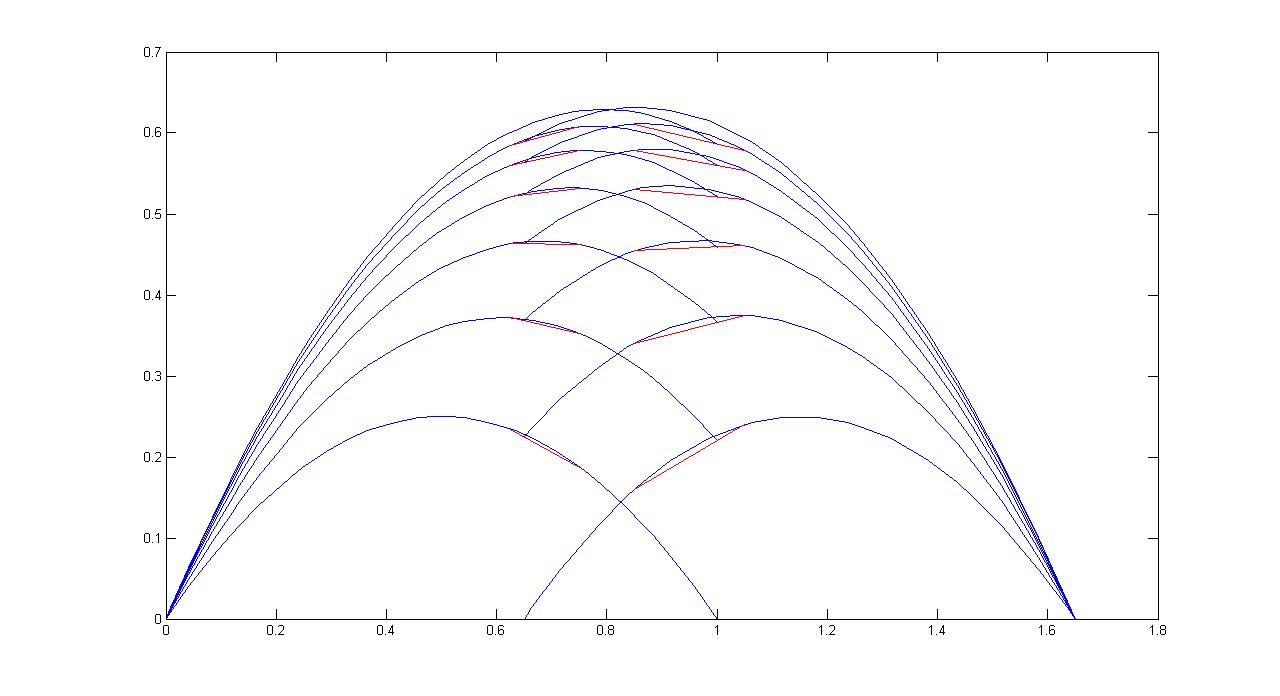}
\caption{One-dimensional Schwarz Method for the equation $u_{xx}=-1$, using the linear interpolation trace operator $P^l$ and zero initial guesses}
\end{center}
\end{figure}

We emphasize here that the matrix 
\begin{equation}
\begin{pmatrix}
\tilde{A}_1 & O \\
O & \tilde{A}_2
\end{pmatrix}
\end{equation}
is not symmetric since the domains may have different number of elements, or even different orders of approximation.

\subsection{Linear interpolation trace operator}
We consider here the linear
interpolation operator we discussed in the previous chapter. 

The algorithm works as was explained above:
consider the first subdomain $\Omega_1 = [0,\beta)$. In order to
construct a linear interpolation of the solution $\widetilde{w}^{n-1}$ at the
point $\eta = \beta$ we need to find the knot span $[\eta_i,
\eta_{i+1})$ which contains $\beta$. We take the values of the
function $\widetilde{w}^{n-1}$ at the ends of this chosen span interval
$\widetilde{w}^{n-1}(\eta_i)$ and $\widetilde{w}^{n-1}(\eta_{i+1})$ and their weighted
sum:
$$
v_{N_{h1}}^{n} = \frac{\beta-\eta_i}{\eta_{i+1} - \eta_{i}}
\widetilde{w}^{n-1}(\eta_i) + \frac{\eta_{i+1}-\beta}{\eta_{i+1} - \eta_{i}}
\widetilde{w}^{n-1}(\eta_{i+1}).
$$

We see that the linear interpolation trace operator is,
actually, a convex sum of the two exact interpolation operators
corresponding to the points $\eta_{i}$ and $\eta_{i+1}$:
\begin{equation}
P_1^l = \frac{\beta-\eta_i}{\eta_{i+1} - \eta_{i}} P^e(\eta_i) +
\frac{\eta_{i+1}-\beta}{\eta_{i+1} - \eta_{i}} P^e(\eta_{i+1}).
\end{equation}
Consequently, the matrix of this interpolation trace operator
can be obtained as the convex sum of the matrices of the exact
trace operators at the points $\eta_{i}$ and $\eta_{i+1}$.

\begin{equation}
\begin{split}
&P_1^l = \frac{\beta-\eta_i}{\eta_{i+1} - \eta_{i}}
\begin{pmatrix}
0 & \dots & 0 & 0 & \dots & 0 & \dots & 0 \\
0 & \dots & 0 & 0 & \dots & 0 & \dots & 0 \\
\vdots & \dots & \vdots & \vdots & \dots & \vdots & \dots & \vdots \\
0 & \dots & \psi_i(\eta_i) & \psi_{i+1}(\eta_i) & \dots &
\psi_{i+p_2}(\eta_i) & \dots & 0
\end{pmatrix}\\
&+\frac{\eta_{i+1}-\beta}{\eta_{i+1} - \eta_{i}}
\begin{pmatrix}
0 & \dots & 0 & 0 & \dots & 0 & \dots & 0 \\
0 & \dots & 0 & 0 & \dots & 0 & \dots & 0 \\
\vdots & \dots & \vdots & \vdots & \dots & \vdots & \dots & \vdots \\
0 & \dots & \psi_{i+1}(\eta_{i+1}) & \psi_{i+2}(\eta_{i+1}) & \dots
& \psi_{i+p_2+1}(\eta_{i+1}) & \dots & 0
\end{pmatrix}.
\end{split}
\end{equation}
The iterative scheme is now:
\begin{equation}
\begin{split}
& \begin{pmatrix}
\tilde{A}_1 & O \\
O & \tilde{A}_2
\end{pmatrix}
\cdot \mathbf{u}^n = f +\begin{pmatrix}
\widetilde{A}_{dir1} & O \\
O & \widetilde{A}_{dir2}
\end{pmatrix}
\cdot 
\begin{pmatrix}
O & P_1^l \\
P_2^l & O
\end{pmatrix}
\cdot \mathbf{u}^{n-1} = \\
&f+\begin{pmatrix}
O & \widetilde{A}_{dir1} P_1^l \\
\widetilde{A}_{dir2} P_2^l & O
\end{pmatrix}
\cdot \mathbf{u}^{n-1},
\end{split}
\end{equation}
where 
\begin{equation}
P_2^l = \frac{\alpha - \xi_j}{\xi_{j+1} - \xi_{j}} P^e(\xi_j) +
\frac{\xi_{j+1}-\alpha}{\xi_{j+1} - \xi_{j}} P^e(\xi_{j+1}).
\end{equation} and $\alpha$ belongs to the knot span $[\xi_j,\xi_{j+1})$.
\subsection{Initial guesses}
Of course, there is some freedom in choosing the initial boundary values $\widetilde{v}^0(\beta)$ and $\widetilde{w}^0(\alpha)$. It turns out that the method converges to the same approximate solutions regardless of the initial guesses since the convergence is guaranteed by the properties of the matrix, which is a M-matrix in our case as shown in the Appendix.

\begin{figure}
\begin{center}
\includegraphics[width=1\textwidth]{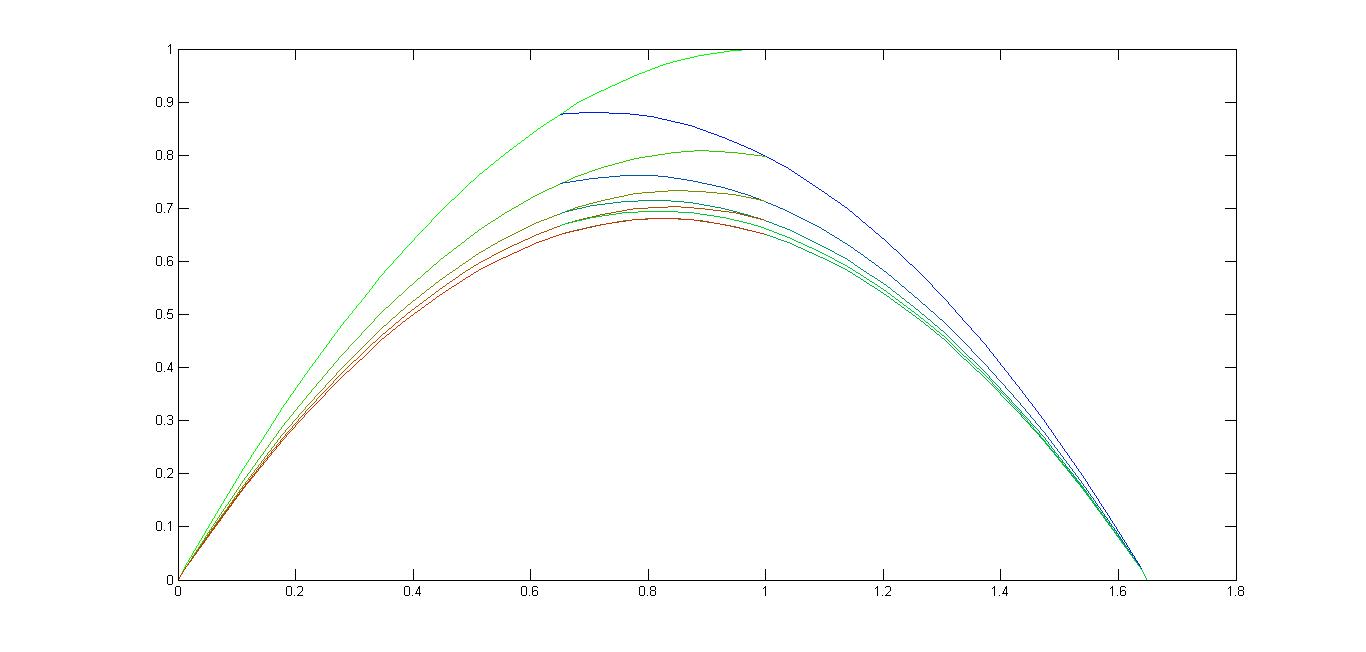}
\caption{Convergence of the one-dimensional SADDM with non-zero initial guesses. Here the initial values for the inner boundary points are taken $\widetilde{v}^0(\beta) = \beta$ and $\widetilde{w}^0(\alpha) = \alpha$}
\label{fig:random_start_iter}
\end{center}
\end{figure}

You can see an example of a convergence procedure for non-zero initial guesses in the Figure \ref{fig:random_start_iter}
\subsection{Multi-dimensional domains}
When we consider a multi - dimensional case we have to perform an approximation of the boundary conditions since the boundary condition function is be a B-spline or NURBS curve or surface in general. This approximation can be performed in one of the ways we have already discussed above in section ~\ref{sec:dir_bc} : quasi-interpolation or least-squares approximation. All the other steps of the solution are similar to the one-dimensional case.
\\
In the next chapter we consider two and three dimensional applications of ASDDM as well as its parallelisation.

\section{Numerical Results}

In this section we provide  illustration of the convergence properties of the ASDDM in our implementation. As we have already mentioned one of the main advantages of the ASDDM methods is that the meshes of different patches are non-matching. This important feature makes it possible to apply this method in CSG and zooming problems. In particular, it allows to perform local zooming at places where the solution possesses weak singularities. Below we demonstrate some examples of this technique and its usage. We also describe how we implement the parallelized version of the solver and show a few three-dimensional multi-patched examples obtained with it. 

As already mentioned before in all the model problems considered in this work the errors are $\mathcal{L}^2$-norm differences between the exact solutions and the approximations, see section (4.1.2).

\section{One-dimensional examples and Analysis}
\subsection{Dependence of the convergence on the B-spline degree}

\begin{figure}
\begin{center}
\includegraphics[width=0.85\textwidth]{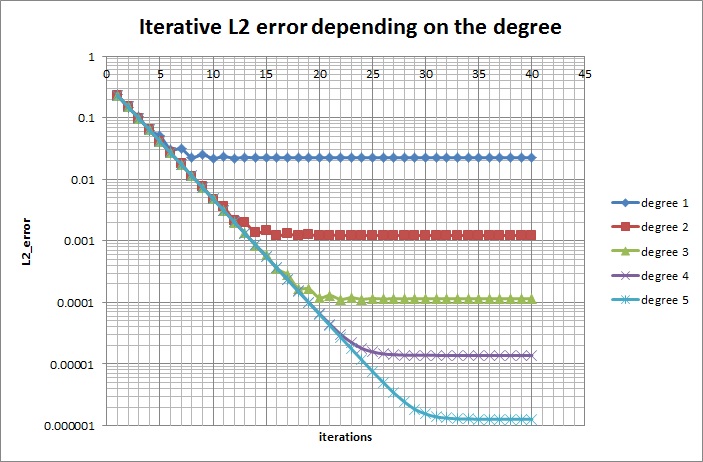}
\label{fig:pdegree}
\caption{$\mathcal{L}^2$-error iterative convergence depending on the polynomial degree, one-dimensional case.}
\end{center}
\end{figure}

First, we treat the one-dimensional problem to show how the iterative process depends on the degree of the B-splines used. The graph given in Figure \ref{fig:pdegree} shows the $\mathcal{L}^2$-errors as a function of the iteration number for different polynomial degrees. Obviously, higher polynomial degrees lead to better precision of the final result, but \emph{the rates of convergence remain the same} which makes the iterative procedure time consuming.

\subsection{Dependence of the convergence on the overlapping area}
\begin{figure}
\begin{center}
\includegraphics[width=0.85\textwidth]{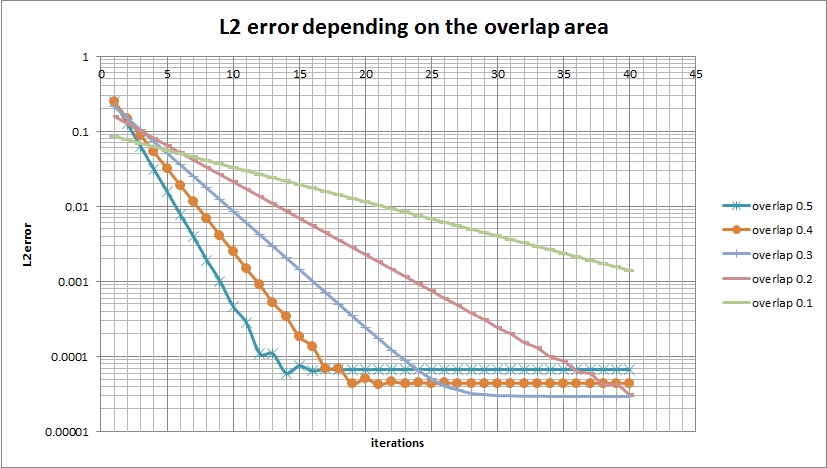}
\caption{$\mathcal{L}^2$-error iterative convergence depending on the ovrelapping area}
\label{fig:overlap}
\end{center}
\end{figure}

Many times when the Domain Decomposition technique is applied to the original domain we have some freedom in choosing the subdomains. For example, when a local zooming is considered the only thing which is predefined is the domain we zoom in, but the choice of the second subdomain is given to us. The question is how the relative overlapping area influences the convergence process. In Figure ~\ref{fig:overlap} we compare different domain decompositions for the same one-dimensional example. The overlapping is measured as the ratio of the overlapping area to the area of one subdomain. Apparently, the larger the overlapping is,  the faster the iterative convergence becomes. An open question is how would one balance the overlapping and the number of degrees of freedom versus the number of iterations to optimize the computation.

\subsection{Dependence of the approximate solution on the mesh size}
\begin{figure}
\begin{center}
\includegraphics[width=0.8\textwidth]{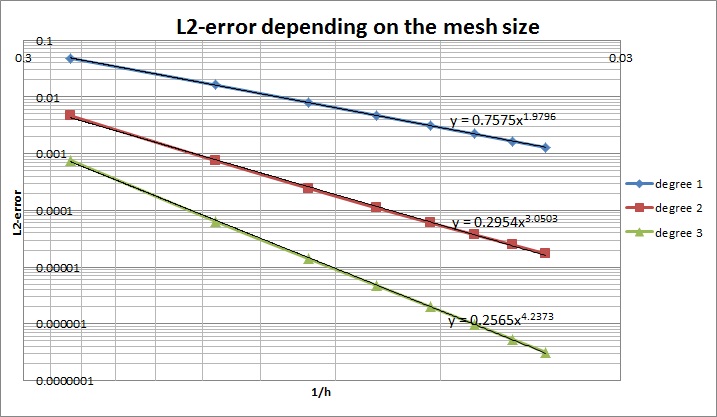}
\caption{One-dimensional $\mathcal{L}^2$-error of the approximate solution depending on the mesh size for different polynomial degrees}
\label{fig:meshsize}
\end{center}
\end{figure}

Figure ~\ref{fig:meshsize} shows how the approximate solution depends on the refinement for different degrees of B-splines. You can see that the dependence is close to theoretically predicted, given by the formula ~\cite{buffa}:

\begin{equation}
\inf_{s \in \mathcal{V}} ||u \circ \mathbf{F} - s||_{\mathcal{L}^2} \leq C h^{p+1} |u \circ \mathbf{F}|_{\mathcal{H}^{p+1}},
\label{conf_formula}
\end{equation}
where $u$ is the exact solution and $C$ is a constant that may depend on the degree $p$ of the polynomials.


\section{Convergence properties in multi-dimensional cases: Application to Zooming}
The formula \ref{conf_formula} becomes different in the multi-dimensional case, but still the approximated solution depends on the mesh size in almost the same way \cite{buffa}. 

Consider the following two-dimensional problem:

Let  $\Omega$ be an open circle of radius 3 placed at the origin:
\begin{equation}
 -\Delta u=\sin(x^2+y^2-9)-4 \cos(x^2+y^2-9) \text{ on } \Omega,
\end{equation}
\begin{equation} 
u|_{\partial\Omega= 0}.
\end{equation}

The exact solution is $u=\sin(x^2+y^2-9)$

\begin{figure}
\begin{center}
\includegraphics[width=0.75\textwidth]{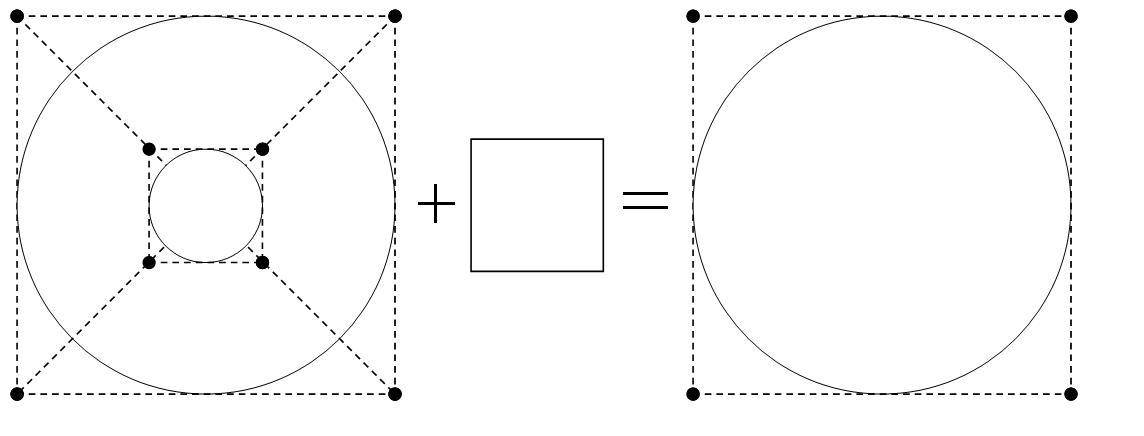}
\caption{Domain Decomposition scheme}
\label{fig:PDD}
\end{center}
\end{figure}

Apply zooming by considering  $\Omega$ as a union of an annulus and a square.
$\Omega = \Omega_{annulus} \cup \Omega_{square}$ as shown in the Figure \ref{fig:PDD}.

\begin{figure}
\begin{center}
\includegraphics[width=0.8\textwidth]{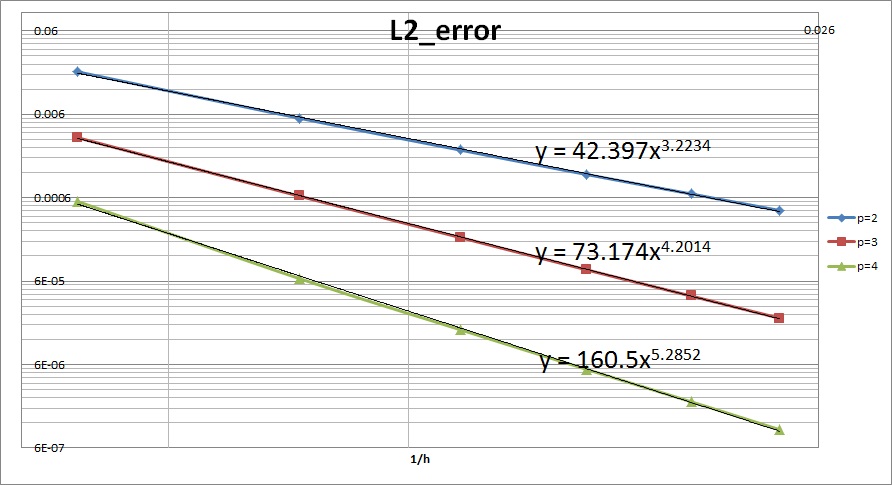}
\caption{Annulus $\mathcal{L}^2$-errors depending on the mesh size for different polynomial degrees}
\label{fig:Pconv}
\end{center}
\end{figure}

\begin{figure}
\begin{center}
\includegraphics[width=0.8\textwidth]{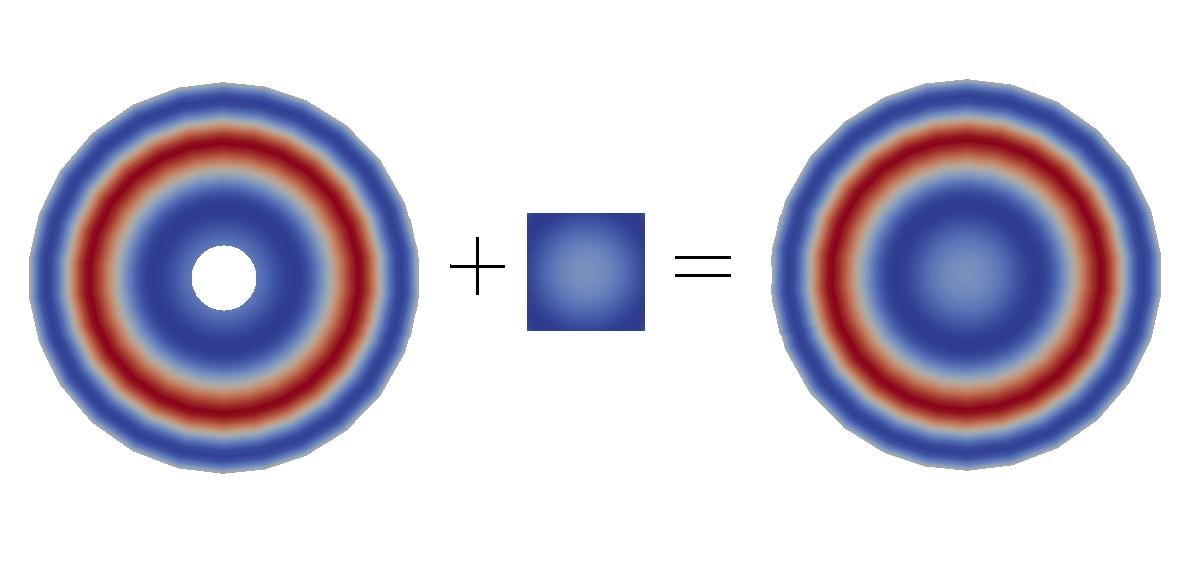}
\caption{Numerical solution}
\label{fig:Psol}
\end{center}
\end{figure}

In Figures \ref{fig:Pconv} and \ref{fig:Psol} below you can see the convergence properties and the numerical solutions we obtained in this case by applying the ASDDM.

\section{Two-dimensional examples}
\subsection{Non-matching domain meshes}
\begin{figure}
\begin{center}
\includegraphics[width=0.75\textwidth]{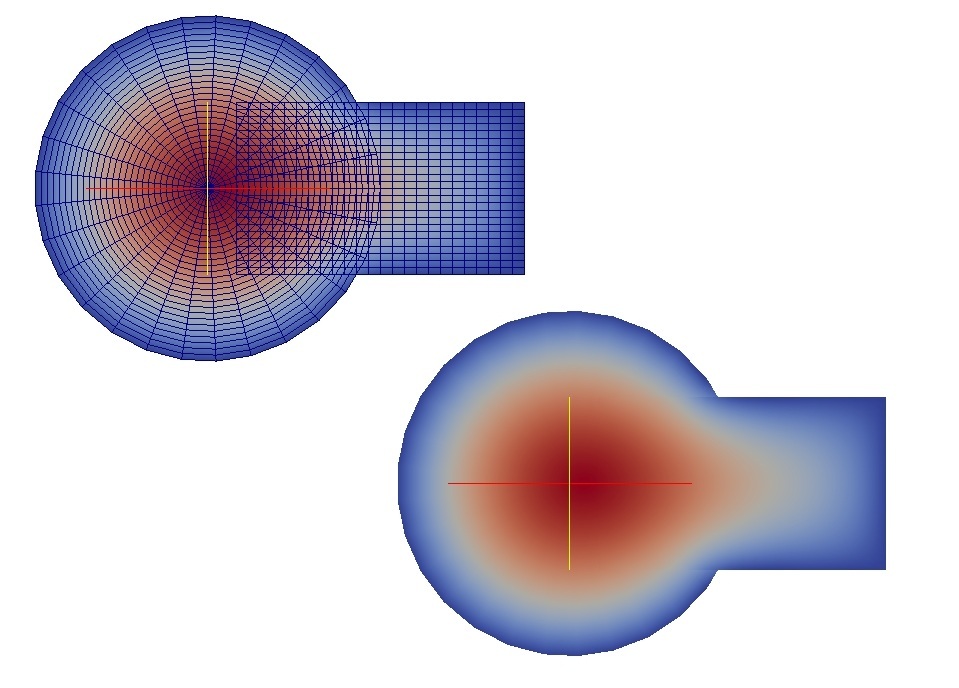}
\caption{Numerical solution over the DD logo.}
\label{fig:logo_numsol}
\end{center}
\end{figure}

An example of non-matching meshes and a solution obtained by the Domain Decomposition method is shown in Figure \ref{fig:logo_numsol}. The domain is taken from the logo of the Domain Decomposition organization, seeFigure ~\ref{fig:logo_s}. We also give in Table 1 the reaction between overlap and DD iterations.
\begin{table}[t]
\caption{Number of  Iteration and max error by Overlap Distance}    
\begin{tabular}{cccccc}
\hline
Overlapp & $0.1$ &$0.25$   &  $ 0.5 $& $0.75$ & $1.0$  \\

Iterations &  max error& at each iteration &  & & \\
\hline
1 & $.6*10^{-1}$  & $.9*10^{-1}$ & $1.2*10^{-1}$& $1.2*10^{-1}$&$ 0.7*10^{-1}$\\
2  & $.3*10^{-1}$  & $.3*10^{-1}$ & $0.16*10^{-1}$& $0.12*10^{-1}$&$ 0.11*10^{-1}$\\
3 & $.15*10^{-1}$  & $.024*10^{-1}$ & $0.01*10^{-1}$& $0.01*10^{-1}$&$ 0.022*10^{-1}$\\
4 & $.05*10^{-1}$  & $.08*10^{-1}$ & $0.034*10^{-1}$& $0.005*10^{-1}$&$ 0.006*10^{-1}$\\
5 & $.03*10^{-1}$  & $.007*10^{-1}$ & $0.006*10^{-1}$& $0.006*10^{-1}$&$ 0.001*10^{-1}$\\
6 & $.02*10^{-1}$  & $.004*10^{-1}$ & $0.002*10^{-1}$& $0.003*10^{-1}$&$ c$\\
7 & $.02*10^{-1}$  & $.002*10^{-1}$ & $0.0015*10^{-1}$& $0.002*10^{-1}$&$ c$\\
8  &$.01*10^{-1}$  & $.002*10^{-1}$ & $0.001*10^{-1}$& $0.001*10^{-1}$&$ c$\\
9&$.006*10^{-1}$  & $.001*10^{-1}$ & $c$& $c $&$ c$\\
10&$.003*10^{-1}$  & $c$ & $c$& $c$&$ c$\\
\hline
\end{tabular}
\label{default}
\end{table}

\subsection{Application of zooming: weak corner singularity}
\begin{figure}
\begin{center}
\includegraphics[width=0.4\textwidth]{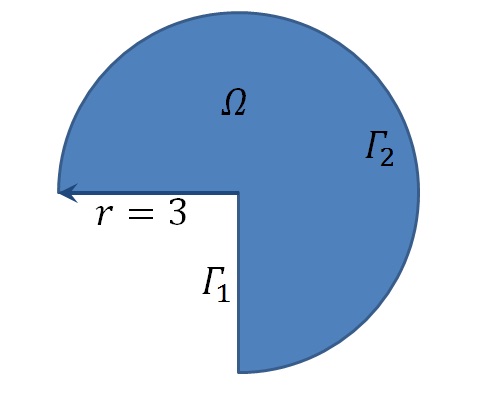}
\caption{Problem with a weak singularity at the corner. The domain.}
\label{fig:givoli_domt}
\end{center}
\end{figure}

We now consider some engineering applications of the Domain Decomposition techniques used for zooming. One of the main applications of ASDDM are the problems where the solution possesses a weak singularity at the corner. The example we treat here is taken from \cite{givoli}.

Given the following domain:
\begin{equation}
\Omega=\{ (\rho,\theta) \in \mathbb{R} | 0 < \rho < 3; -\frac{\pi}{2} < \theta < \pi \},
\end{equation}
which is an open circular sector of angle $\alpha = \frac{3\pi}{2}$ as shown in  Figure \ref{fig:givoli_domt} and the equation:

\begin{equation}
\Delta u = 0 \text{ in } \Omega,
\end{equation}

\begin{displaymath}
   u = \left\{
     \begin{array}{ll}
       0 & : \text{on }\Gamma_1 \subset \Gamma_\Omega = \partial{\Omega}\\
       \theta (\alpha - \theta) & : \text{on }\Gamma_2  = \Gamma_\Omega \backslash \Gamma_1.
     \end{array}
   \right.
\end{displaymath} 

\begin{figure}
\begin{center}
\includegraphics[width=0.8\textwidth]{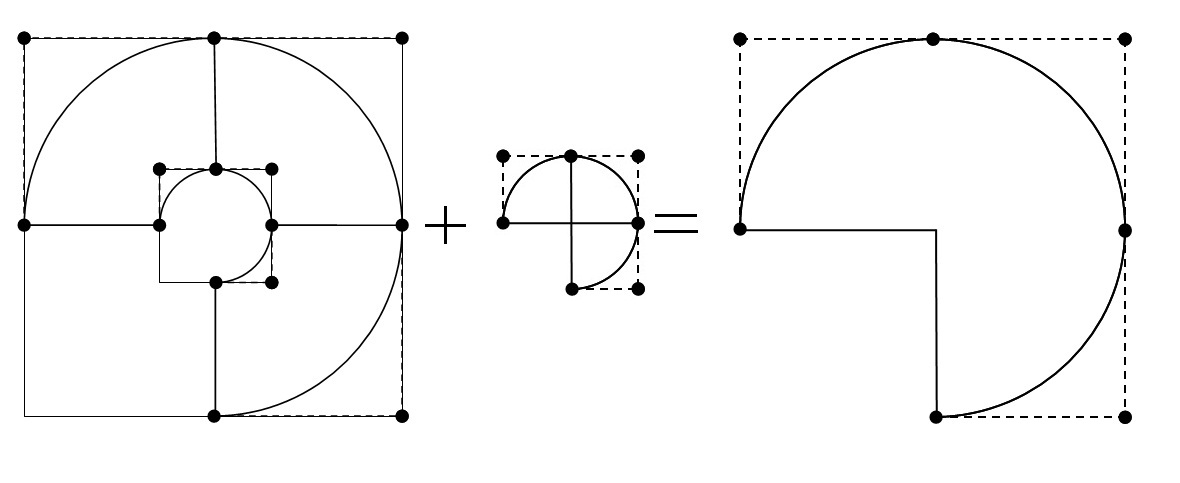}
\caption{NURBS definition, Boundary and Domain Decomposition}
\label{fig:givoli_tog_1}
\end{center}
\end{figure}

The exact solution is given by the following uniformly convergent series:
\begin{equation}
u_{ex}(\rho, \theta) = \frac{9}{\pi} \sum_{n=1,3,5,\dots} \frac{1}{n^3} \left (\frac{\rho}{r}\right)^\frac{2 n}{3}\sin \left(\frac{2 n \theta}{3} \right). 
\label{eq:ex_sol_giv}
\end{equation}
We see from ~\eqref{eq:ex_sol_giv} that the radial derivative is unbounded when approaching the origin which may cause the numerical solution to be unstable at the vicinity of the origin. In order to zoom the solution near the corner we decompose the domain in the following way: $\Omega = \Omega_{out} \cup \Omega_{zoom},$ as shown in Figure \ref{fig:givoli_tog_1}.

\begin{figure}
\begin{center}
\includegraphics[width=0.75\textwidth]{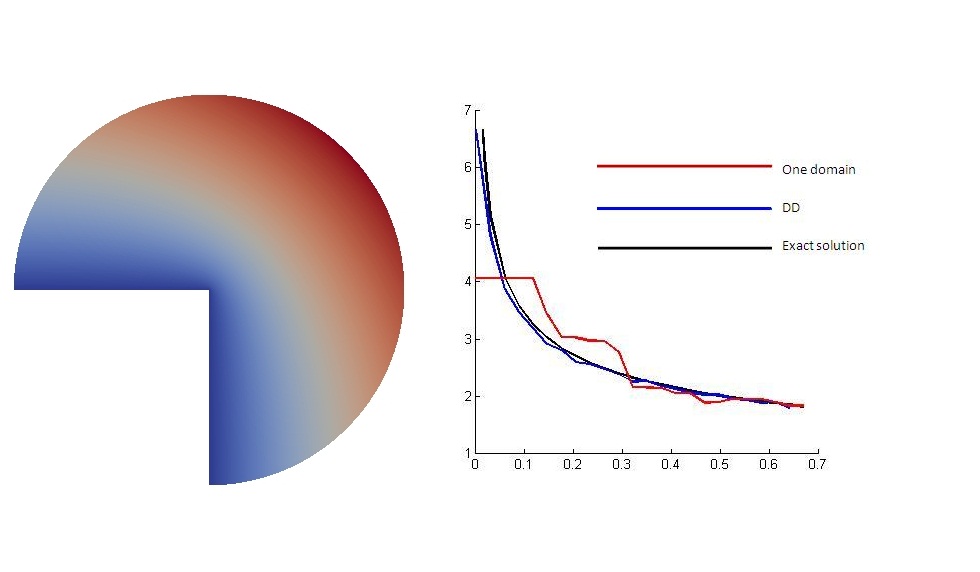}
\caption{Numerical solution of (3.2) and the corresponding  radial derivative}
\label{fig:givoli_derivative}
\end{center}
\end{figure}

Figure \ref{fig:givoli_derivative} shows the numerical solution and the radial derivative calculated for this problem. The principal order of the radial derivative of the exact solution is $-\frac{1}{3}$, and the numerically computed result using the Domain Decomposition technique coincides with this value.

\begin{figure}
\begin{center}
\includegraphics[width=0.7\textwidth]{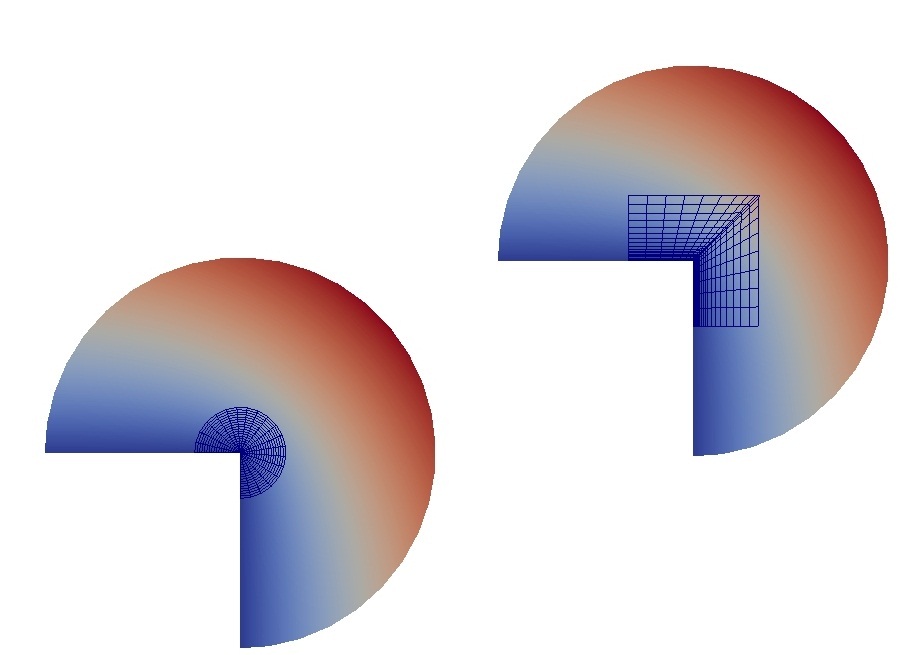}
\caption{Different zooming regions}
\label{fig:givoli_different}
\end{center}
\end{figure}

We have tested different forms of the zooming region as shown in Figure \ref{fig:givoli_different} and in all the cases the results coincided.

\begin{figure}
\begin{center}
\includegraphics[width=0.8\textwidth]{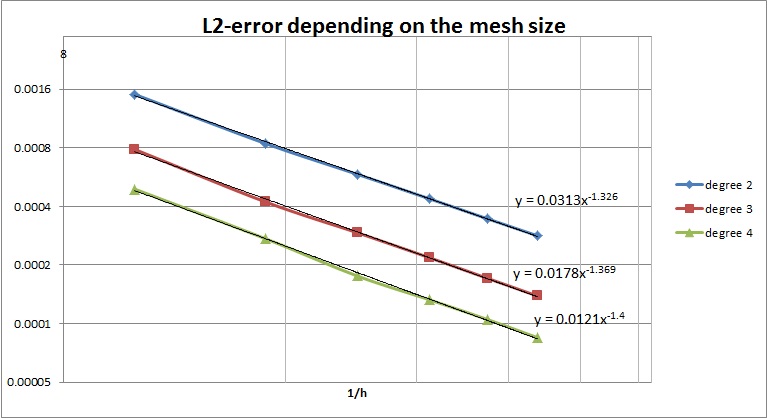}
\caption{$\mathcal{L}^2$-error of the approximate solution depending on the mesh size for different polynomial degrees}
\label{fig:givoli_conv}
\end{center}
\end{figure}

The solution in this case is not regular and possesses a weak singularity at the corner and the error estimate \ref{conf_formula} does not apply anymore. The dependence of the approximate solution on the mesh size becomes different from the regular case and the rate of convergence does not improve with greater polynomial degree as shown in Figure ~\ref{fig:givoli_conv}.

\section{Three-dimensional examples}
\begin{figure}
\begin{center}
\includegraphics[width=0.7\textwidth]{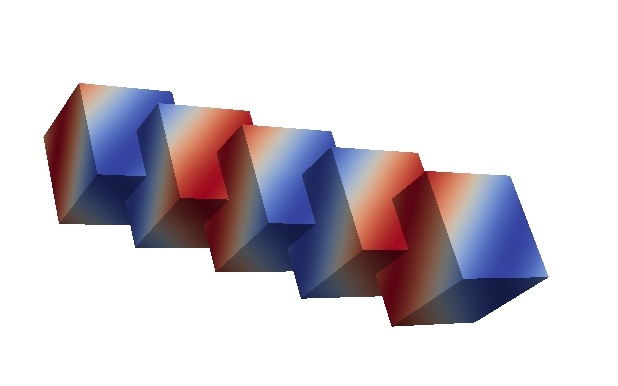}
\caption{Chain of cubes with the exact solution $\sin(x+y+z)$}
\label{fig:cubes}
\end{center}
\end{figure}

\begin{figure}
\begin{center}
\includegraphics[width=0.7\textwidth]{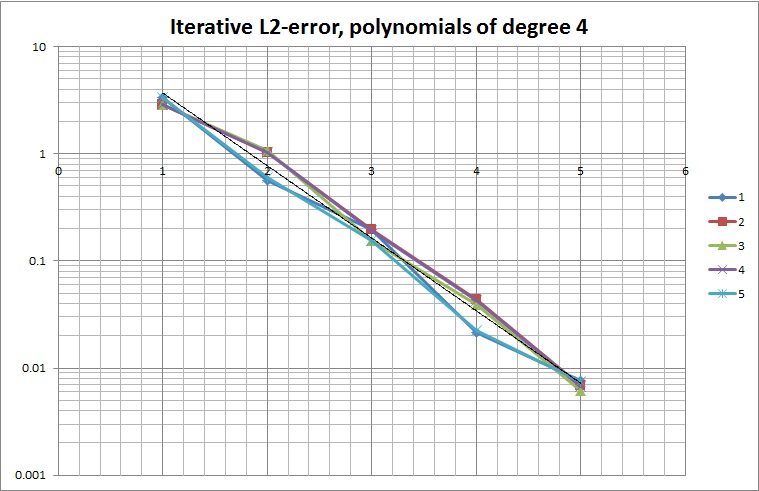}
\caption{Iterative $\mathcal{L}^2$-error for the chain of cubes}
\label{fig:cubes_conv}
\end{center}
\end{figure}

In this section we present some examples obtained for multi-patched bodies.
Consider the equation:
\begin{equation}
 -\Delta u=\sin(x+y+z) \text{ on } \Omega,
\end{equation}
\begin{equation} 
u|_{\partial\Omega}=\sin(x+y+z),
\end{equation}
where $\Omega$ is a chain of 5 overlapping cubes, as shown in Figure \ref{fig:cubes}. The exact solution is $u=\sin(x+y+z)$.
You can see the rates of convergence of the solution for all the 5 cubes in the Figure \ref{fig:cubes_conv}.

\subsection{Parallelized algorithm and multi-patched problems}
As we have mentioned in the beginning of this work, one of the main advantages of the ASDDM is that it may be easily parallelized. The modern computers are all multi-core which provides us with an opportunity of real parallelization of the code, meaning that the solution on each patch may be computed independently and simultaneously on a different processor. 

The algorithm we implemented works as following. Before we start the iterations we precompute all the date that the solver needs, then we open a worker thread for each subprocess, corresponding to a patch, and start the iterative procedure. On each iteration the $n$-th thread solves the equation on the $n$-th subdomain according to the scheme shown in Figure \ref{fig:algo_struct} of the previous chapter. After the $n$-th iteration computation is finished, the data is synchronized between the threads. This synchronization is necessary to impose the correct boundary conditions for each subdomain on the next iteration.

\begin{figure}
\begin{center}
\includegraphics[width=0.8\textwidth]{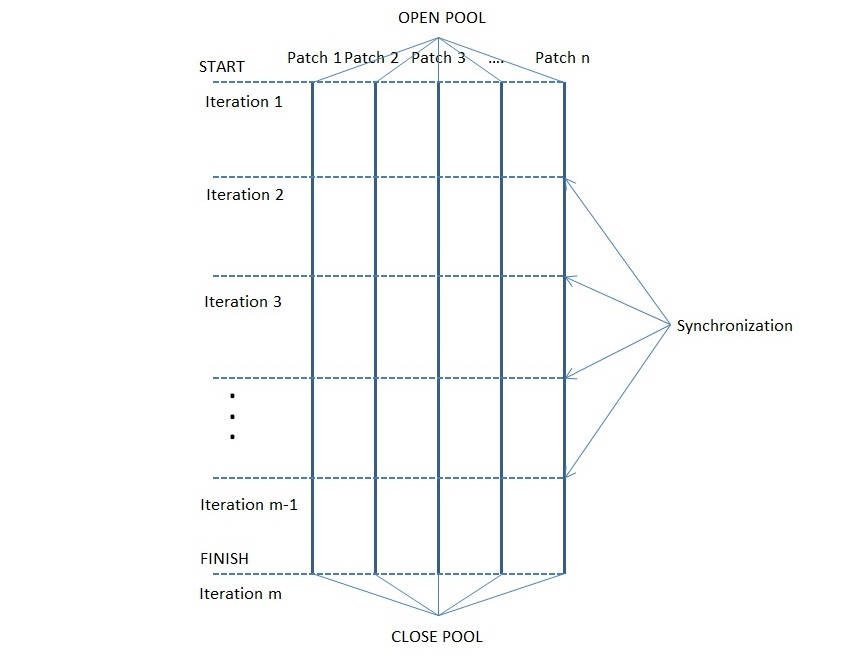}
\caption{Synchronization algorithm}
\label{fig:syncro}
\end{center}
\end{figure}

The schematic synchronization mechanism is shown in Figure \ref{fig:syncro}.

 \subsection{Elasticity example for multi-patched domain}
\begin{figure}
\begin{center}
\includegraphics[width=0.9\textwidth]{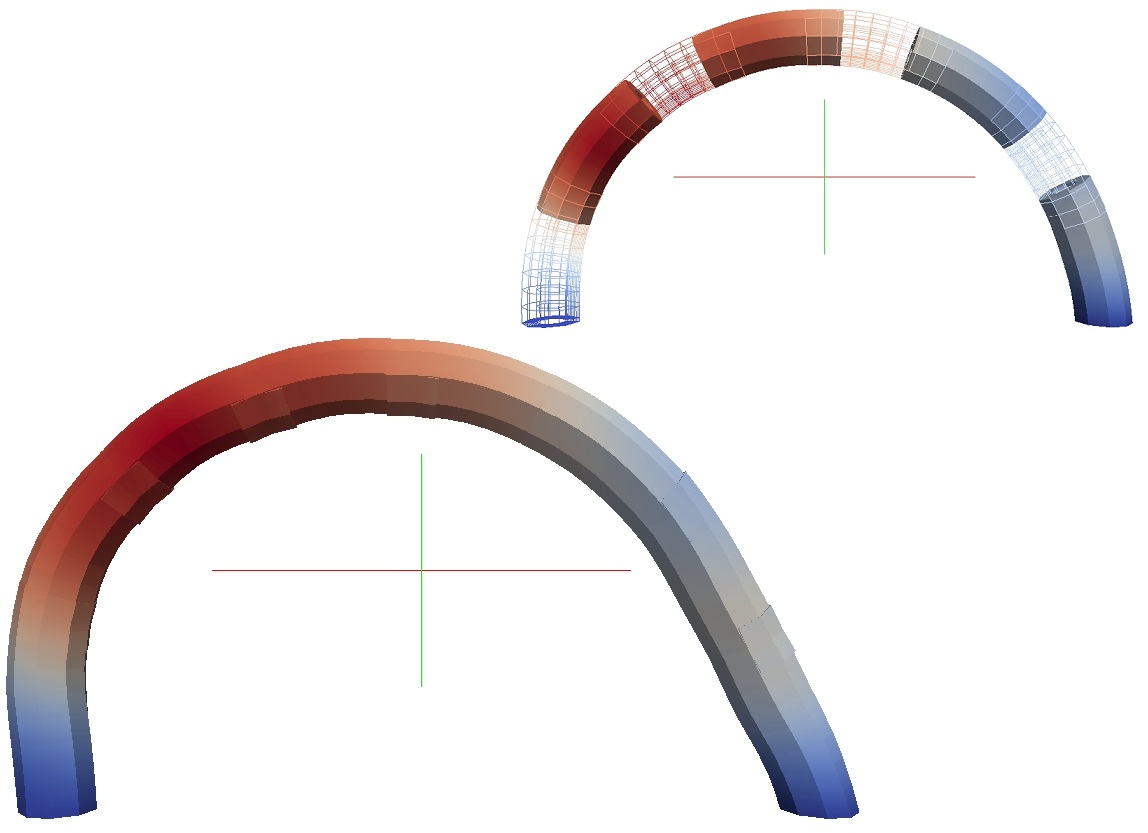}
\caption{Elasticity example}
\label{fig:pipes}
\end{center}
\end{figure}

As an example of parallel computation we took an elasticity problem. The domain is a hollow thick half-ring with fixed ends which is placed into a gravitational field. This patch was decomposed into 8 overlapping subdomains as shown in the Figure \ref{fig:pipes} and the solution was obtained using the parallelized method described above.

We performed numerical tests on the 4-processors machine and achieved a full 4-times enhancement in the time of computations.

\subsection{Statistics} In this section we illustrate the dependence of the number of degrees of freedom in the problem as a function of dimension, refinement and polynomial degrees. Of course, these numbers do not depend on the isogeometrical mapping, so we calculated the data using a cube domain. 

In the two-dimensional case we obtained:
\vspace{.15in}
\begin{center}
\begin{tabular}{cccc}
\hline\noalign{\smallskip}
\multicolumn{4}{r}{degree of the B-splines} \\
\cline{2-4}\noalign{\smallskip}
Partitions of the mesh & $2$ & $3$ & $4$ \\
\noalign{\smallskip}\hline\noalign{\smallskip}
0  & 10  & 27 & 52\\
1  & 27  & 52 & 85\\
2  & 52  & 85 & 126\\
5  & 175 & 232 & 297\\
10 & 540 & 637 & 742\\
11 & 637 & 742 & 855\\
20 & 1870 & 2047 & 2232\\
\noalign{\smallskip}\hline
\end{tabular}
\end{center}
\vspace{.3in}

In the three-dimensional case we obtained:

\vspace{.15in}
\begin{center}
\begin{tabular}{ccc}
\hline\noalign{\smallskip}
\multicolumn{3}{r}{degree of the B-splines} \\
\cline{2-3}\noalign{\smallskip}
Partitions of the mesh & $1$ & $2$ \\
\noalign{\smallskip}\hline\noalign{\smallskip}
1  & 765 & 624\\
2  & 2100 & 1275\\
\noalign{\smallskip}\hline
\end{tabular}
\end{center}

The number of degrees of freedom, and consequently, the computational time grows rapidly with the dimension of the problem, especially if the mesh need to be refined. Here the Domain Decomposition techniques become very helpful.

\section{Conclusion} The main purpose of this work was to establish the connection and application of the recently discovered IGA framework to solve PDEs using Additive Domain Decomposition Method. It turns out that IGA may be successfully used for solving different problems even on complex geometries, thus giving us all the advantages of IGA without restrictions on the geometry. Another important application of these methods are zooming problems. The IGA framework is restricted to its native techniques like T-splines and here the Domain Decomposition becomes very useful. 

More research is needed is to explore  ways of approximating and imposing the Dirichlet boundary conditions. Different techniques of imposing boundary conditions may influence the computational times and the precision of the approximate solution significantly hence a research subject is improving techniques to define such BCs on trimmed patches.

The question of optimal overlap choice (  a large overlap seems better, but more costly in degrees of freedom) needs more study. Also what is the best way to define the solution on the intersections is an opened question. 

Another important direction for  research is the preconditioning of the solver in large scale problems. An efficient pre-conditioner may improve the computational times and simplify the structure of the solver( and insure convergence , since we do not have a maximum principle here.)

Finally a CSG construct is not made only of unions of primitives, and it my be the results of Boolean differences, substractions etc.,thus DD should be coupled with other methods such as fictitious domains .

\bibliography{bib1}
\bibliographystyle{plain}

\appendix

\let\clearpage\relax
\eject

\end{document}